\documentclass[final,leqno,onefignum,onetabnum]{scrartcl}
\usepackage[T1]{fontenc}
\usepackage[latin9]{inputenc}
\usepackage{verbatim}
\usepackage{float}
\usepackage{amsmath}
\usepackage{amssymb}
\usepackage{graphicx}
\usepackage{esint}

\makeatletter

\floatstyle{ruled}
\newfloat{algorithm}{tbp}{loa}
\providecommand{\algorithmname}{Algorithm}
\floatname{algorithm}{\protect\algorithmname}

\usepackage{texdraw}
\usepackage{placeins}
\usepackage{epsfig}
\@ifundefined{definecolor}
 {\usepackage{color}}{}
\usepackage{listings}
\usepackage{tikz}
\usepackage{pgf}
\usetikzlibrary{arrows}

\newtheorem{prop}{Proposition}[section]

\makeatother

\title{Image Segmentation with Eigenfunctions of an Anisotropic Diffusion Operator
\thanks{This work was supported in part by the NSF under Grant DMS-1115118.}
} 

\author{
Jingyue Wang\thanks{Department of Mathematics, 
the University of Kansas, Lawrence, KS 66045, U.S.A. 
({\em jwang@math.ku.edu}).}
\and 
Weizhang Huang\thanks{Department of Mathematics, 
the University of Kansas, Lawrence, KS 66045, U.S.A. 
({\em whuang@ku.edu}).}
}

\begin{document}
\maketitle
\newcommand{\slugmaster}{%
\slugger{siims}{xxxx}{xx}{x}{x--x}}

\begin{abstract}
We propose the eigenvalue problem of an anisotropic diffusion operator for image segmentation.
The diffusion matrix is defined based on the input image. The eigenfunctions and the projection
of the input image in some eigenspace capture key features of the input image.
An important property of the model is that for many input images, the first few eigenfunctions are close to being piecewise
constant, which makes them useful as the basis for a variety of applications such as image segmentation and edge detection.
The eigenvalue problem is shown to be related to the algebraic eigenvalue problems resulting from
several commonly used discrete spectral clustering models.
The relation provides a better understanding and helps developing more efficient numerical implementation
and rigorous numerical analysis for discrete spectral segmentation methods.
The new continuous model is also different from energy-minimization methods such as geodesic active contour
in that no initial guess is required for in the current model. The multi-scale feature is a natural consequence of
the anisotropic diffusion operator so there is no need to solve the eigenvalue problem at multiple levels.
A numerical implementation based on a finite element method with an anisotropic mesh adaptation strategy
is presented. It is shown that the numerical scheme gives much more accurate results on eigenfunctions than
uniform meshes. Several interesting features of the model are examined in numerical examples and
possible applications are discussed.
\end{abstract}

\noindent{\bf Key Words.}
eigenvalue problem, image segmentation, finite-element schemes, mesh adaptation, anisotropic diffusion.

\noindent{\bf AMS 2010 Mathematics Subject Classification.}
65N25, 68U10, 94A08

\pagestyle{myheadings}
\thispagestyle{plain}
\markboth{IMAGE SEGMENTATION WITH EIGENFUNCTIONS}{IMAGE SEGMENTATION WITH EIGENFUNCTIONS}

\section{Introduction}

We are concerned with image segmentation using
the eigenvalue problem of an anisotropic linear diffusion operator,
\begin{equation}
-\nabla \cdot (\mathbb{D}\nabla u)=\lambda u, \quad \text{ in }\Omega
\label{eq:HW-eigen}
\end{equation}
subject to a homogeneous Dirichlet or Neumann boundary condition, where the diffusion matrix
$\mathbb{D}$ is symmetric and uniformly positive definite on $\Omega$. In this study, we consider
an anisotropic diffusion situation where $\mathbb{D}$ has different eigenvalues and is defined based
on the gray level of the input image. 


A method employing an eigenvalue problem to study image segmentation is referred to as a spectral clustering
method in the literature. This type of methods have extracted great interest from researchers
in the past decade; e.g., see \cite{Grady-01,shi_normalized_2000,Wu-Leahy-1993}.
They are typically derived from a minimum-cut criterion  on a graph.
One of the most noticeable spectral clustering methods is the normalized cut method proposed
by Shi and Malik \cite{shi_normalized_2000} (also see Section~\ref{SEC:relation-1} below)
which is based on the eigenvalue problem
\begin{equation}
(D-W) {\bf u} = \lambda D {\bf u} ,
\label{eq:Malik-Shi-eigen}
\end{equation}
where ${\bf u}$ is a vector representing the gray level value on the pixels, $W$ is a matrix defining
pairwise similarity between pixels, and $D$ is a diagonal matrix formed with the degree
of pixels (cf. Section 2.2 below). The operator $L=D-W$ corresponds to the graph Laplacian in graph spectral
theory. An eigenvector associated with the second eigenvalue is used
as a continuous approximation to a binary or $k$-way vector that
indicates the partitions of the input image. Shi and Malik suggested that
image segmentation be done on a hierarchical basis where
low level coherence of brightness, texture, and etc. guides a binary (or
$k$-way) segmentation that provides a big picture while high
level knowledge is used to further partition the low-level segments.


While discrete spectral clustering methods give impressive partitioning results in general,
they have several drawbacks. Those methods are typically defined and operated
on a graph or a data set. Their implementation cost depends on the size of the graph or data set.
For a large data set, they can be very expensive to implement.
Moreover, since they are discrete, sometimes their physical and/or geometrical meanings
are not so clear. As we shall see in Section~\ref{SEC:relation-1}, the normalized cut
method of Shi and Malik \cite{shi_normalized_2000} is linked to an anisotropic diffusion
differential operator which from time to time can lead to isotropic diffusion.

The objective of this paper is to investigate the use of the eigenvalue problem
(\ref{eq:HW-eigen}) of an anisotropic diffusion operator for image segmentation.
This anisotropic model can be viewed as a continuous, improved anisotropic generalization
of discrete spectral clustering models such as (\ref{eq:Malik-Shi-eigen}).
The model is also closely related to the Perona-Malik anisotropic filter.
The advantages of using a continuous model for image segmentation
include (i) It has a clear physical interpretation (heat diffusion or
Fick's laws of diffusion in our case);
(ii) Many well developed theories of partial differential equations can be used;
(iii) Standard discretization methods such as finite differences, finite elements,
finite volumes, and spectral methods can be employed;
and (iv) The model does not have be discretized on a mesh associated with the given
data set and indeed, mesh adaptation can be used to improve accuracy and efficiency.
As mentioned early, we shall define the diffusion matrix $\mathbb{D}$ using
the input image and explore properties of the eigenvalue problem.
One interesting property is that for many input images, the first few eigenfunctions of the model
are close to being piecewise constant, which are very useful for image segmentation.
However, this also means that these eigenfunctions change abruptly between
objects and their efficient numerical approximation requires mesh adaptation.
In this work, we shall use an anisotropic mesh adaptation strategy developed
by the authors \cite{Huang-Wang-13} for differential eigenvalue problems.
Another property of (\ref{eq:HW-eigen}) is that 
eigenfunctions associated with small eigenvalues possess coarse, global
features of the input image whereas eigenfunctions associated with larger eigenvalues
carry more detailed, localized features.
The decomposition of features agrees with the view of Shi and Malik \cite{shi_normalized_2000}
on the hierarchical structure of image segmentation but in a slightly different sense
since all eigenfunctions come from low level brightness knowledge.

The paper is organized as follows. In Section~\ref{SEC:eigen}, we give a detailed description
of the eigenvalue problem based on an anisotropic diffusion operator and
discuss its relations to some commonly used discrete spectral clustering models 
and diffusion filters and some other models in image segmentation.
Section~\ref{SEC:implement} is devoted to 
the description of the finite element implementation of the model and
an anisotropic mesh adaptation strategy.
In Section~\ref{SEC:numerics}, we present a number of applications in image segmentation
and edge detection and demonstrate several properties of the model.
Some explanations to the piecewise constant property of eigenfunctions are given in
Section~\ref{SEC:piecewise}.
Concluding remarks are given in Section~\ref{SEC:conclusion}.

\section{Description of the eigenvalue problem}
\label{SEC:eigen}

\subsection{Eigenvalue problem of an anisotropic diffusion operator}

We shall use the eigenvalue problem (\ref{eq:HW-eigen}) subject to a Dirichlet or Neumann boundary condition
for image segmentation.


We are inspired by the physics of anisotropic heat transport
process (e.g., see \cite{Gunter-Yu-Kruger-Lackner-05,Sharma-Hammett-07}),
treating the dynamics of image diffusion as the transport of energy
(pixel values) and viewing the eigenvalue problem as the steady state of the dynamic process.
Denote the principal diffusion direction by
$v$ (a unit direction field) and its perpendicular unit direction by $v^{\perp}$.
Let the conductivity coefficients along these directions be
$\chi_{\parallel}$ and $\chi_{\perp}$. ($v$, $\chi_{\parallel}$, and $\chi_{\perp}$ will be defined below.) 
Then the diffusion matrix can be written as
\begin{equation}
\mathbb{D}=\chi_{\parallel}vv^{T}+\chi_{\perp}v^{\perp}(v^{\perp})^{T} .
\label{D-1}
\end{equation}
When $\chi_{\parallel}$ and $\chi_{\perp}$ do not depend on $u$, the diffusion operator
in (\ref{eq:HW-eigen}) is simply a linear symmetric second order elliptic operator.
The anisotropy of the diffusion tensor $\mathbb{D}$ depends on the choice of the conductivity coefficients.
For example, if $\chi_{\parallel}\gg\chi_{\perp}$, the diffusion
is preferred along the direction of $v$. Moreover, if $\chi_{\parallel}=\chi_{\perp}$,
the diffusion is isotropic, having no preferred diffusion direction.

To define $\mathbb{D}$, we assume that an input image is given. Denote its gray level by $u_0$.
In image segmentation, pixels with similar values of gray level will be grouped and
the interfaces between those groups provide object boundaries. Since those interfaces are orthogonal 
to $\nabla u_0$, it is natural to choose the principal diffusion direction as $v=\nabla u_{0}/|\nabla u_{0}|$.
With this choice, we can rewrite (\ref{D-1}) into
\begin{equation}
\mathbb{D} = \frac{\chi_{\parallel}}{|\nabla u_{0}|^{2}}
\begin{bmatrix}|\partial_x u_{0}|^{2}+\mu |\partial_y u_{0}|^{2} & (1-\mu)\left|\partial_x u_{0}\partial_y u_{0}\right|\\
(1-\mu)\left|\partial_x u_{0} \partial_y u_{0}\right| & |\partial_y u_{0}|^{2}+\mu |\partial_x u_{0}|^{2} \end{bmatrix}
\label{D-2}
\end{equation}
where $\mu = \chi_{\perp}/\chi_{\parallel}$.
We consider two choices of $\chi_{\parallel}$ and $\mu$. The first one is
\begin{equation}
\chi_{\parallel}=g(|\nabla u_{0}|), \quad \mu = 1, 
\label{D-3}
\end{equation}
where $g(x)$ is a conductance function that governs the behavior of diffusion.
This corresponds to linear isotropic diffusion. As in \cite{Perona-Malik-90},
we require $g$ to satisfy $g(0)=1$, $g(x)\ge0$, and $g(x)\to 0$ as $x \to \infty$.
For this choice, both $\chi_{\parallel}$ and $\chi_{\perp}$ becomes very small
across the interfaces of the pixel groups and therefore, almost no diffusion is allowed
along the normal and tangential directions of the interfaces. 

The second choice is 
\begin{equation}
\chi_{\parallel}=g(|\nabla u_{0}|), \quad \mu=1+|\nabla u_{0}|^{2}.
\label{D-4}
\end{equation}
This choice results in an anisotropic diffusion process. Like the first case, almost no diffusion is allowed
across the interfaces of the pixel groups but, depending on the choice of $g$, some degree of diffusion
is allowed on the tangential direction of the interfaces.
We shall show later that with a properly chosen $g$ the eigenfunctions of (\ref{eq:HW-eigen})
capture certain ``grouping'' features of the input image $u_{0}$
very well. This phenomenon has already been observed and explored
in many applications such as shape analysis \cite{Reuter-09,Reuter-06},
image segmentation and data clustering \cite{Grady-01, shi_normalized_2000, Shi-Malik-2001, Wu-Leahy-1993},
and high dimensional data analysis and machine learning
\cite{Belkin_towards_2005,Nadler_diffusion_2005,Nadler_diffusion_2006,Luxburg-2007}.
In these applications, all eigenvalue problems are formulated on a
discrete graph using the graph spectral theory, which is different from what is
considered here, i.e., eigenvalue problems of differential operators.
The application of the latter to image segmentation is much less known.
We shall discuss the connection of these discrete eigenvalue problems
with continuous ones in the next subsection.

It is noted that the gray level function $u_0$ is defined only at pixels.
Even we can view $u_{0}$ as the ``ground truth'' function (assuming
there is one function whose discrete sample is the input image),
it may not be smooth and the gradient cannot be defined
in the classical sense. Following \cite{Alvarez-Lions-92,Catte-Lions-92},
we may treat $u_{0}$ as a properly regularized approximation of the ``true image''
so that the solution to the eigenvalue problem (\ref{eq:HW-eigen}) exists.
In the following, we simply take $u_{0}$ as the linear interpolation
of the sampled pixel values (essentially an implicit regularization
from the numerical scheme). More sophisticated regularization
methods can also be employed. 

We only deal with gray level images in this work. The approach
can be extended to color or texture images when a diffusion matrix
can be defined appropriately based on all channels. In our computation, we use
both Dirichlet and Neumann boundary conditions, with the latter
being more common in image processing. 

\subsection{Relation to discrete spectral clustering models}
\label{SEC:relation-1}

The eigenvalue problem (\ref{eq:HW-eigen}) is closely related to a family of discrete spectral
clustering models, with the earliest one being the normalized
cut method proposed by Shi and Malik \cite{shi_normalized_2000}.
To describe it, we define the degree of dissimilarity (called $cut$) between
any two disjoint sets $A,B$ of a weighted undirected graph $G=(V,E)$ (where $V$ and $E$ denote the sets of the nodes
and edges of the graph) as the total weight of
the edges connecting nodes in the two sets, i.e., 
\[
cut(A,B)=\sum_{p\in A,\; q\in B}w(p,q).
\]
Wu and Leahy \cite{Wu-Leahy-1993} proposed to find $k$-subgraphs by
minimizing the maximum cut across the subgroups and use them for
a segmentation of an image. However, this approach
usually favors small sets of isolated nodes in the graph. To address
this problem, Shi and Malik \cite{shi_normalized_2000} used
the normalized cut defined as 
\[
Ncut(A,B)=\frac{cut(A,B)}{assoc(A,A\cup B)}+\frac{cut(A,B)}{assoc(B,A\cup B)} ,
\]
where $assoc(A,A\cup B)=\sum_{p\in A,\; q\in A\cup B}w(p,q)$. They sought the minimum
of the functional $Ncut(A,B)$ recursively to obtain a $k$-partition
of the image. The edge weight $w(p,q)$ is chosen as
\[
w(p,q)=\begin{cases}
e^{-|u_{q}-u_{p}|^2/\sigma^2}, & q\in\mathcal{N}_{p},\\
0, & {\rm otherwise,}
\end{cases}
\]
where $\mathcal{N}_{p}$ is a neighborhood of pixel $p$ and $\sigma$ is a positive parameter.
Shi and Malik showed that the above optimization problem is NP-hard
but a binary solution to the normalized cut problem can be
mapped to a binary solution to the algebraic eigenvalue problem (\ref{eq:Malik-Shi-eigen})
with $D$ being a diagonal matrix with diagonal entries $d_{p}=\sum_{q}w(p,q)$
and $W$ being the weight matrix $(w(p,q))_{p,q}^{N\times N}$. Eigenvectors of
this algebraic eigenvalue problem are generally not binary. They are used to approximate
binary solutions of the normalized cut problem through certain partitioning.

To see the connection between the algebraic eigenvalue problem (\ref{eq:Malik-Shi-eigen})
(and therefore, the normalized cut method) with the continuous eigenvalue problem (\ref{eq:HW-eigen}),
we consider an eigenvalue problem in the form of (\ref{eq:HW-eigen}) with the diffusion matrix defined as
\begin{equation}
\mathbb{D} = \begin{bmatrix}e^{-|\partial_{x}u_{0}|^2/\sigma^2} & 0\\
0 & e^{-|\partial_{y}u_{0}|^2/\sigma^2} \end{bmatrix}
\label{D-5}
\end{equation}
A standard central finite difference discretization of this problem on a rectangular mesh gives rise to
\begin{equation}
\frac{(c_{E_{i,j}}+c_{W_{i,j}}+c_{N_{i,j}}+c_{S_{i,j}})u_{i,j}-c_{E_{i,j}}u_{i+1,j}-c_{W_{i,j}}u_{i-1,j}-c_{N_{i,j}}u_{i,j+1}-c_{S_{i,j}}u_{i,j-1}}{h^{2}}=\lambda u_{i,j},
\label{eq:orthotropic-fe-scheme}
\end{equation}
where $h$ is the grid spacing and the coefficients $c_{E_{i,j}},\; c_{W_{i,j}},\; c_{N_{i,j}},\; c_{S_{i,j}}$ are given as
\begin{eqnarray*}
c_{E_{i,j}} = e^{-|u_{i+1,j}-u_{i,j}|^2/\sigma^2},\quad
c_{W_{i,j}} = e^{-|u_{i-1,j}-u_{i,j}|^2/\sigma^2},\\
c_{N_{i,j}} = e^{-|u_{i,j+1}-u_{i,j}|^2/\sigma^2},\quad
c_{S_{i,j}} =  e^{-|u_{i,j-1}-u_{i,j}|^2/\sigma^2}.
\end{eqnarray*}
It is easy to see that (\ref{eq:orthotropic-fe-scheme}) is almost the same as (\ref{eq:Malik-Shi-eigen})
with the neighborhood $\mathcal{N}_{i,j}$ of a pixel location $(i,j)$ being chosen to include the four closest
pixel locations $\{(i+1,j),(i-1,j),(i,j+1),(i,j-1)\}$. The difference lies in that (\ref{eq:Malik-Shi-eigen}) has
a weight function on its right-hand side. Moreover, it can be shown that (\ref{eq:orthotropic-fe-scheme})
gives {\it exactly} the algebraic eigenvalue problem for the average cut problem
\[
{\rm min}\frac{cut(A,B)}{|A|}+\frac{cut(A,B)}{|B|} ,
\]
where $|A|$ and $|B|$ denote the total numbers of nodes in $A$ and $B$, respectively.
Notice that this problem is slightly different from the normalized cut problem and its solution is known
as the Fiedler value. Furthermore, if we consider the following generalized eigenvalue problem
(by multiplying the right-hand side of (\ref{eq:HW-eigen}) with a mass-density function),
\begin{equation}
-\nabla \cdot \left(\begin{bmatrix}e^{-|\partial_{x}u_{0}|^2/\sigma^2} & 0\\
0 & e^{-|\partial_{y}u_{0}|^2/\sigma^2} \end{bmatrix}
\nabla u\right)= (e^{-|\partial_x u_{0}|^2/\sigma^2} + e^{-|\partial_y u_{0}|^2/\sigma^2}) \lambda u, 
\label{eq:pm-aniso-1}
\end{equation}
we can obtain (\ref{eq:Malik-Shi-eigen}) exactly with a proper central finite difference discretization.

The above analysis shows that either the average cut or normalized cut model can be approximated
by a finite difference discretization of the continuous eigenvalue problem
(\ref{eq:HW-eigen}) with the diffusion matrix (\ref{D-5})
which treats diffusion differently in the $x$ and $y$ directions.
While (\ref{D-5}) is anisotropic in general, it results in isotropic diffusion
near oblique interfaces where $\partial_x u_0 \approx \partial_y u_0$ or 
$\partial_x u_0 \approx - \partial_y u_0$. This can be avoided
with the diffusion matrix (\ref{D-2}) which defines diffusion differently along
the normal and tangential directions of group interfaces.
In this sense, our method consisting of (\ref{eq:HW-eigen}) with (\ref{D-2})
can be regarded as an improved version of (\ref{eq:HW-eigen}) with (\ref{D-5}),
and thus, an improved continuous generalization of the normalized cut or the average cut method.

It should be pointed out that there is a fundamental difference between discrete spectral clustering
methods and those based on continuous eigenvalue problems.
The former are defined and operated directly on a graph or data set
and their cost depends very much on the size of the graph or data.
On the other hand, methods based on continuous eigenvalue problems
treat an image as a sampled function and are defined by
a discretization of some differential operators. They have the advantage
that many standard discretization methods such as finite difference, finite
element, finite volume, and spectral methods can be used. 
Another advantage is that they do not have to be operated directly
on the graph or the data set. As shown in \cite{Huang-Wang-13},
continuous eigenvalue problems can be solved efficiently on
adaptive, and especially anisotropic adaptive, meshes (also see Section~\ref{SEC:numerics}).

It is worth pointing out that the graph Laplacian
can be connected to a continuous diffusion operator by defining the latter on a
manifold and proving it to be the limit of the discrete Laplacian.
The interested reader is referred to the work of
\cite{Belkin_towards_2005,Nadler_diffusion_2005,Nadler_diffusion_2006,Singer-06,Luxburg-2007}.


\subsection{Relation to diffusion models}

The eigenvalue problem (\ref{eq:HW-eigen}) is related to several diffusion models used in image processing.
They can be cast in the form
\begin{equation}
\frac{\partial u}{\partial t}=\nabla \cdot \left(\mathbb{D}\nabla u\right) 
\label{eq:linear-diffusion}
\end{equation}
with various definitions of the diffusion matrix. For example, the Perona-Malik nonlinear filter
\cite{Perona-Malik-90} is in this form with $\mathbb{D} = g(|\nabla u|) I$, where $g$ is the same
function in (\ref{D-3}) and $I$ is the identity matrix. The above equation with $\mathbb{D}$ defined in (\ref{D-2})
with $\mu=1$ and $\chi_{\parallel}=g(|\nabla u_{0}|)$ gives rise to a linear diffusion process
that has similar effects as the affine Gaussian smoothing process \cite{Nitzberg-Shiota-92}.
The diffusion matrix we use in this paper in most cases is in the form (\ref{D-2}) with $\mu$ and $\chi_{\parallel}$
defined in (\ref{D-4}). A similar but not equivalent process
was studied as a structure adaptive filter by Yang et al. \cite{Yang-Burger-96}.
The diffusion matrix (\ref{D-2}) can be made $u$-dependent by choosing $\mu$ and $\chi_{\parallel}$
as functions of $\nabla u$.
Weickert \cite{Weickert-1996} considered a nonlinear anisotropic diffusion model with a diffusion matrix
in a similar form as (\ref{D-2}) but with $\nabla u_0$ being replaced by the gradient of a smoothed
gray level function $u_\sigma$ and with  $\chi_{\parallel} = g(|\nabla u_\sigma|)$
and $\mu = 1/g(|\nabla u_\sigma|)$. 

Interestingly, Perona and Malik \cite{Perona-Malik-90} considered 
\begin{eqnarray}
\frac{\partial u}{\partial t}  = \nabla \cdot \left(\begin{bmatrix}g(|\partial_{x}u |) & 0\\
0 & g(|\partial_{y} u|) \end{bmatrix}\nabla u\right) 
\label{eq:pm-linear-diffusion}
\end{eqnarray}
as an easy-to-compute variant to the Perona-Malik diffusion model (with $\mathbb{D} = g(|\nabla u|) I$).
Zhang and Hancock in \cite{Zhang-Hancock-08} considered 
\begin{eqnarray}
\frac{\partial u}{\partial t} = -\mathcal{L}(u_0) u,
\label{eq:Zhang-Hancock}
\end{eqnarray}
where $\mathcal{L}$ is the graph Laplacian defined on the input image $u_0$
and image pixels are treated as the nodes of a graph. The weight between
two nodes $i,j$ is defined as
\[
w_{i,j}=\begin{cases}
e^{-(u_{0}(i)-u_{0}(j))^{2}/\sigma^{2}},  & \quad \text{ for }\|i-j\|\le r\\
0, & \quad \text{otherwise}
\end{cases}
\]
where $r$ is a prescribed positive integer and $\sigma$ is a positive parameter.
As in Section~\ref{SEC:relation-1}, it can be shown that this model can be regarded
as a discrete form of a linear anisotropic diffusion model.
It has been reported in \cite{Buades-Chien-08, Nitzberg-Shiota-92, Yang-Burger-96,Zhang-Hancock-08}
that the image denoising effect with this type of linear diffusion model is
comparable to or in some cases better than nonlinear evolution models.

\section{Numerical Implementation}
\label{SEC:implement}

The eigenvalue problem (\ref{eq:HW-eigen}) is discretized using the standard
linear finite element method with a triangular mesh for $\Omega$. The finite element
method preserves the symmetry of the underlying continuous problem
and can readily be implemented with (anisotropic) mesh adaptation.
As will be seen in Section~\ref{SEC:numerics}, the eigenfunctions
of (\ref{eq:HW-eigen}) can have very strong anisotropic behavior,
and (anisotropic) mesh adaptation is essential to improving the efficiency
of their numerical approximation.

While both Dirichlet and Neumann boundary conditions are considered in our
computation, to be specific we consider only a Dirichlet boundary condition
in the following. The case with a Neumann boundary condition can be discussed similarly.

We assume that a triangular mesh $\mathcal{T}_{h}$ is given for $\Omega$.
Denote the number of the elements of $\mathcal{T}_h$ by $N$ and
the linear finite element space associated with $\mathcal{T}_h$ by $V^{h}\subset H_{0}^{1}\left(\Omega\right)$.
Then the finite element approximation to the eigenvalue problem (\ref{eq:HW-eigen}) subject
to a Dirichlet boundary condition is to find $0 \not\equiv u^h \in V^{h}$ and $\lambda^h \in \mathbb{R}$ such that
\begin{equation}
\int_{\Omega}(\nabla v^h)^t \mathbb{D}\nabla u^h  =\lambda^h \int_{\Omega} u^h v^h,\qquad\forall v^h\in V^{h}.
\label{eq:fem-1}
\end{equation}
This equation can be written into a matrix form as
\[
A {\bf u} = \lambda^h M {\bf u},
\]
where $A$ and $M$ are the stiffness and mass matrices, respectively, and ${\bf u}$ is the vector
formed by the nodal values of the eigenfunction at the interior mesh nodes.

An error bound for the linear finite element approximation of the eigenvalues is given by
a classical result of Raviart and Thomas \cite{Raviart-Thomas-83}. It states that
for any given integer $k$ ($1\le k\le N$),
\[
0\le\frac{\lambda_{j}^{h}-\lambda_{j}}{\lambda_{j}^{h}}\le C(k)\sup_{v\in E_{k},\|v\|=1}\| v-\Pi_{h}v\|_E^{2},
\qquad1\le j\le k
\]
where $\lambda_j$ and $\lambda_j^h$ are the eigenvalues (ordered in an increasing order) of the continuous and
discrete problems, respectively, $E_{k}$ is the linear space spanned by the first $k$ eigenfunctions of the
continuous problem, $\Pi_{h}$ is the projection operator from $L^{2}(\Omega)$ to the finite element space $V^{h}$,
and $\| \cdot \|_E$ is the energy norm, namely,
\[
\| v-\Pi_{h}v\|_E^{2} = \int_\Omega \nabla (v-\Pi_{h}v)^t \mathbb{D} \nabla (v-\Pi_{h}v).
\]
It is easy to show (e.g., see \cite{Huang-Wang-13}) that the project error can be bounded
by the error of the interpolation associated with the underlying finite element space,
with the latter depending directly on the mesh. When the eigenfunctions change abruptly
over the domain and exhibit strong anisotropic behavior, anisotropic mesh adaptation
is necessary to reduce the error or improve the computational efficiency (e.g. see \cite{Boff-10,Huang-Wang-13}).

An anisotropic mesh adaptation method was proposed for eigenvalue problems
by the authors \cite{Huang-Wang-13}, following the so-called $\mathbb{M}$-uniform mesh approach
developed in \cite{Huang-05,Huang-06,Huang-Russell-11} for the numerical solution
of PDEs. Anisotropic mesh adaptation provides one advantage over isotropic one
in that, in addition to the size, the orientation of triangles is also adapted to be aligned with
the geometry of the solution locally. In the context of image processing, this mesh alignment
will help better capture the geometry of edges than with isotropic meshes.
The $\mathbb{M}$-uniform mesh approach of anisotropic mesh adaptation views
and generates anisotropic adaptive meshes as uniform ones in the metric specified
by a metric tensor $\mathbb{M} = \mathbb{M}(x,y)$.
Putting it in a simplified scenario, we may consider a uniform mesh defined on the surface of
the gray level $u$ and obtain an anisotropic adaptive mesh by projecting the uniform mesh
into the physical domain. In the actual computation, instead of using the surface of $u$ we employ
a manifold associated with a metric tensor defined based on the Hessian of the eigenfunctions.
An optimal choice of the metric tensor (corresponding to the energy norm) is given \cite{Huang-Wang-13} as
\[
\mathbb{M}_{K}=\det\left(H_{K}\right)^{-1/4}\max_{(x,y)\in K}\|H_{K}\mathbb{D}(x,y)\|^{1/2}\left(\frac{1}
{|K|}\|H_{K}^{-1}H\|_{L^{2}(K)}^{2}\right)^{1/2}H_{K},\qquad\forall K\in\mathcal{T}_{h}
\]
where $K$ denotes a triangle element of the mesh, $H$ is the intersection of the recovered Hessian
matrices of the computed first $k$ eigenfunctions, and $H_{K}$ is the average of $H$ over $K$. 
A least squares fitting method is used for Hessian recovery. That is, 
a quadratic polynomial is constructed locally for each node via least squares  fitting to neighboring
nodal function values and then an approximate Hessian at the node is obtained by differentiating the
polynomial. The recovered Hessian is regularized with a prescribed small positive constant which is taken
to be $0.01$ in our computation.

An outline of the computational procedure of the anisotropic adaptive mesh finite
element approximation for the eigenvalue problem (\ref{eq:HW-eigen}) is given
in Algorithm~\ref{alg:aniso}. In Step 5, BAMG (Bidimensional Anisotropic Mesh Generator)
developed by Hecht \cite{Hecht-Bamg-98} is used to generate the new mesh
based on the computed metric tensor defined on the current mesh.
The resultant algebraic eigenvalue problems are solved using the Matlab eigenvalue solver
{\tt eigs} for large sparse matrices. Note that the algorithm is iterative. Ten iterations are used in our
computation, which was found to be enough to produce an adaptive mesh
with good quality (see \cite{Huang-05} for mesh quality measures).

\begin{algorithm}[h]
\begin{raggedright}
1. Initialize a background mesh.
\par\end{raggedright}

\begin{raggedright}
2. Compute the stiffness and mass matrices on the mesh.
\par\end{raggedright}

\begin{raggedright}
3. Solve the algebraic eigenvalue problem for the first $k$ eigenpairs.
\par\end{raggedright}

\begin{raggedright}
4. Use the eigenvectors obtained in Step 3 to compute the metric tensor.
\par\end{raggedright}

\begin{raggedright}
5. Use the metric tensor to generate a new mesh (anisotropic, adaptive) and go to Step 2.
\par\end{raggedright}

\caption{Anisotropic adaptive mesh finite element approximation for eigenvalue problems.}
\label{alg:aniso} 
\end{algorithm}

\section{Numerical results}
\label{SEC:numerics}

In this section, all input images are of size $256\times256$ and
normalized so that the gray values are between 0 and 1. The domain
of input images is set to be $[0,1]\times[0,1]$. All eigenfunctions
are computed with a homogeneous Neumann boundary condition unless otherwise specified.
When we count the indices of eigenfunctions, we ignore the first trivial
constant eigenfunction and start the indexing from the second one.

\subsection{Properties of eigenfunctions}

\subsubsection{Almost piecewise constant eigenfunctions}

\begin{figure}
\begin{centering}
\includegraphics[width=10cm]{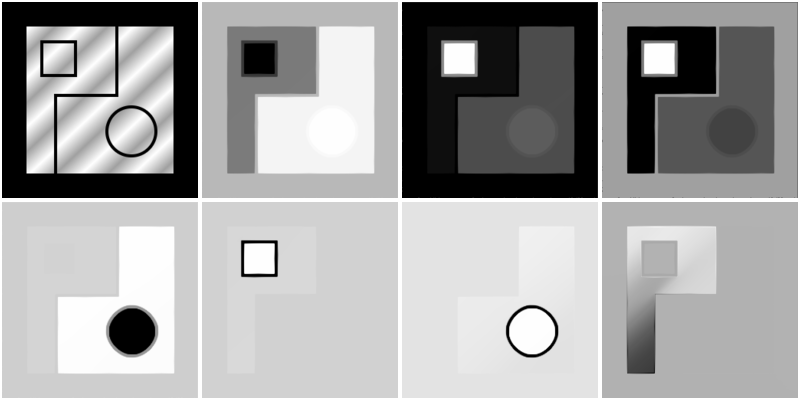}
\par\end{centering}
\caption{The input synthetic image and the first 7 eigenfunctions (excluding
the trivial constant eigenfunction), from left to right, top to bottom. The results are obtained with
the diffusion matrix defined in (\ref{D-2}), (\ref{D-4}), and (\ref{D-6}) ($\alpha = 1.5$).}
\label{fig:synth1}
\end{figure}

A remarkable feature of the eigenvalue problem (\ref{eq:HW-eigen}) with the diffusion matrix (\ref{D-2})
is that for certain input images, the first few eigenfunctions are close
to being piecewise constant. In Fig.~\ref{fig:synth1}, we display
a synthetic image containing 4 objects and the first 7 eigenfunctions.
The gaps between objects are 4 pixel wide. To make the problem more interesting,
the gray level is made to vary within each object (so the gray value of the 
input image is not piecewise-constant).
We use the anisotropic diffusion tensor $\mathbb{D}$ defined in (\ref{D-2}) and (\ref{D-4}) with
\begin{equation}
g(x)=\frac{1}{(1+x^{2})^{\alpha}}, 
\label{D-6}
\end{equation}
where $\alpha$ is a positive parameter. Through numerical experiment (cf. Section~\ref{SEC:4.1.6}), we observe that
the larger $\alpha$ is, the closer to being piecewise constant the eigenfunctions are.
In the same time, the eigenvalue problem (\ref{eq:HW-eigen}) is also harder to solve numerically since
the eigenfunctions change more abruptly between the objects. 
We use $\alpha = 1.5$ in the computation for Fig.~\ref{fig:synth1}.
The computed eigenfunctions are normalized such that 
they have the range of $[0,255]$ and can be rendered as gray level images. 
The results are obtained with an adaptive mesh of 65902 vertices and re-interpolated to
a $256\times256$ mesh for rendering.

\begin{figure} 
\begin{minipage}{.2\linewidth} 
\includegraphics[scale=.5]{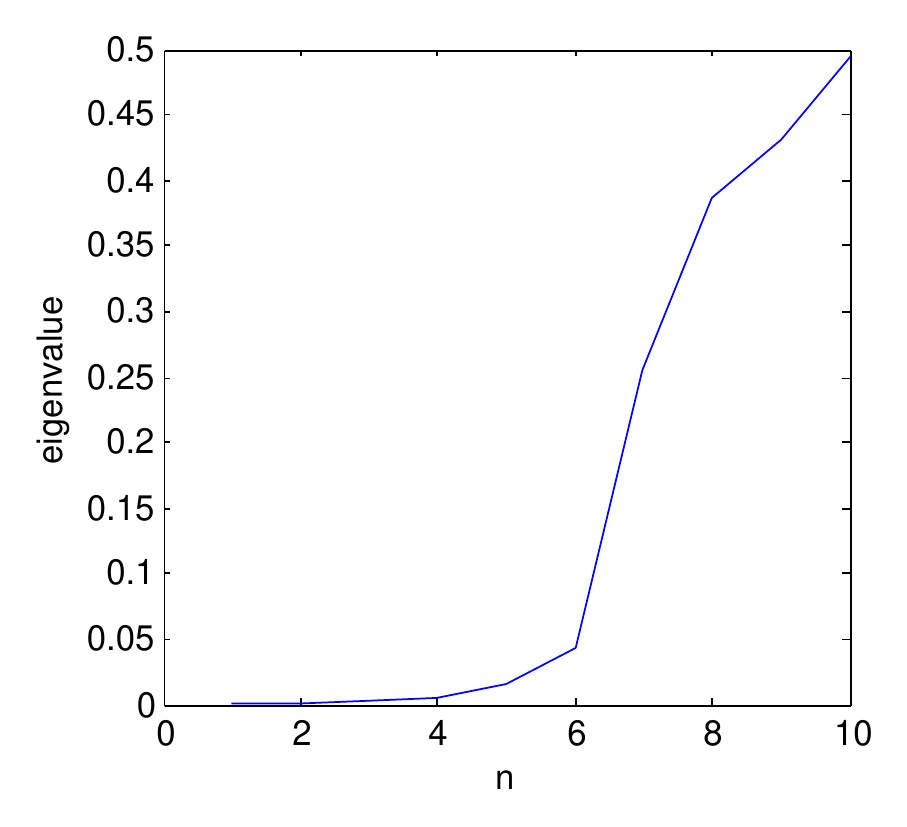}
\end{minipage}%
\hfill
\begin{minipage}{.2\linewidth} 
\centering 
\includegraphics[scale=.23]{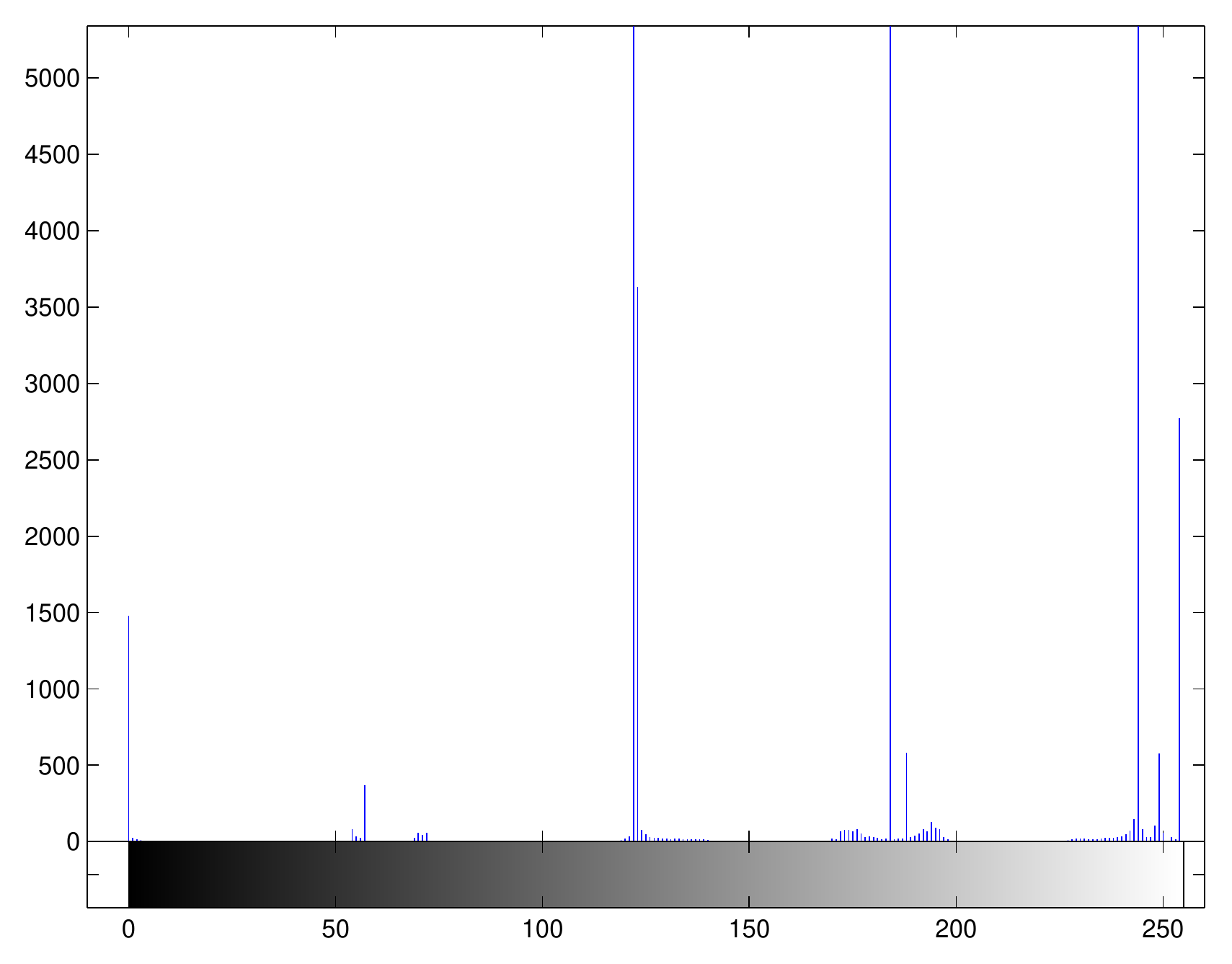} 
\end{minipage} 
\hfill
\begin{minipage}{.2\linewidth} 
\centering 
\includegraphics[scale=.23]{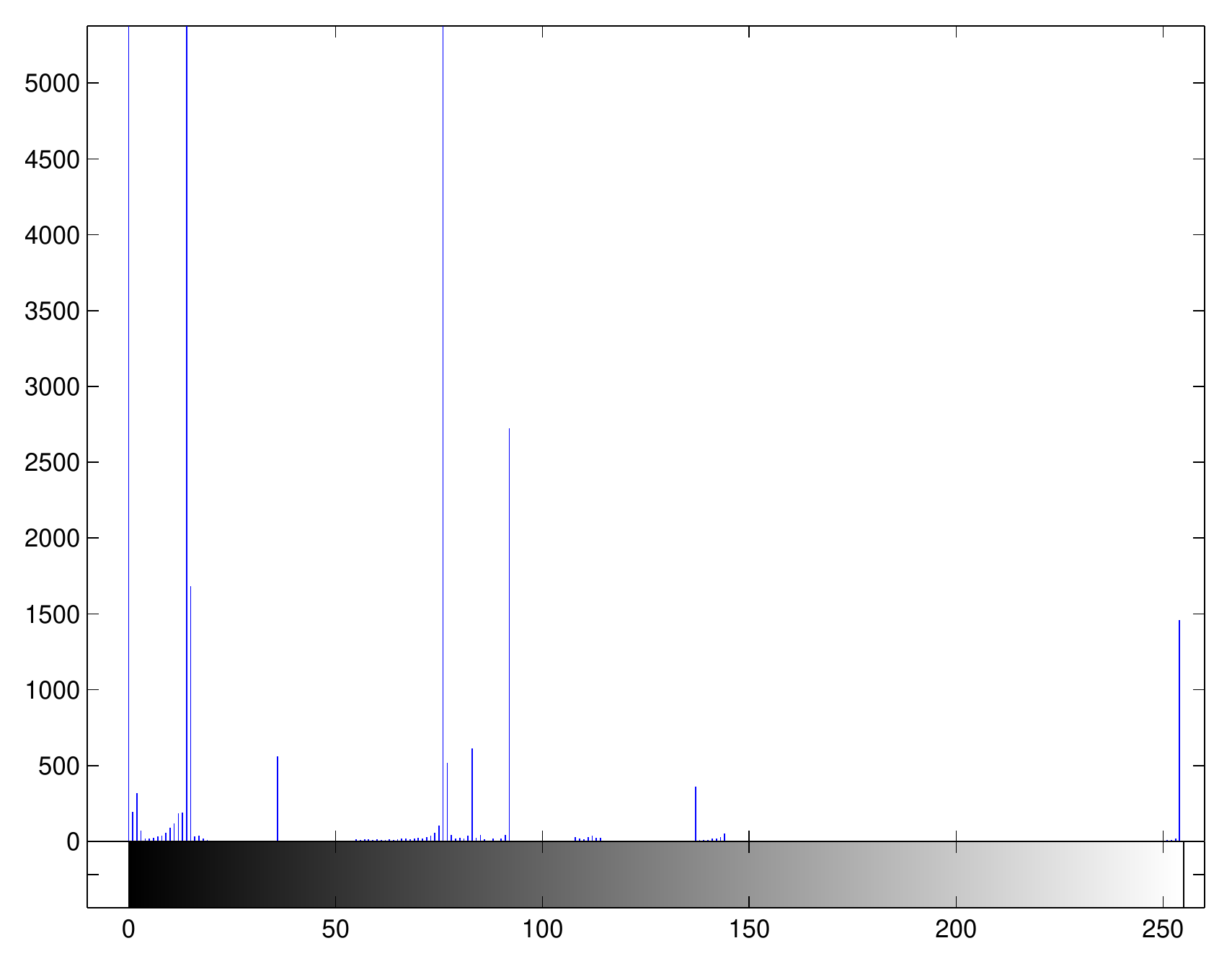} 
\end{minipage} 
\hfill
\begin{minipage}{.2\linewidth} 
\centering 
\includegraphics[scale=.23]{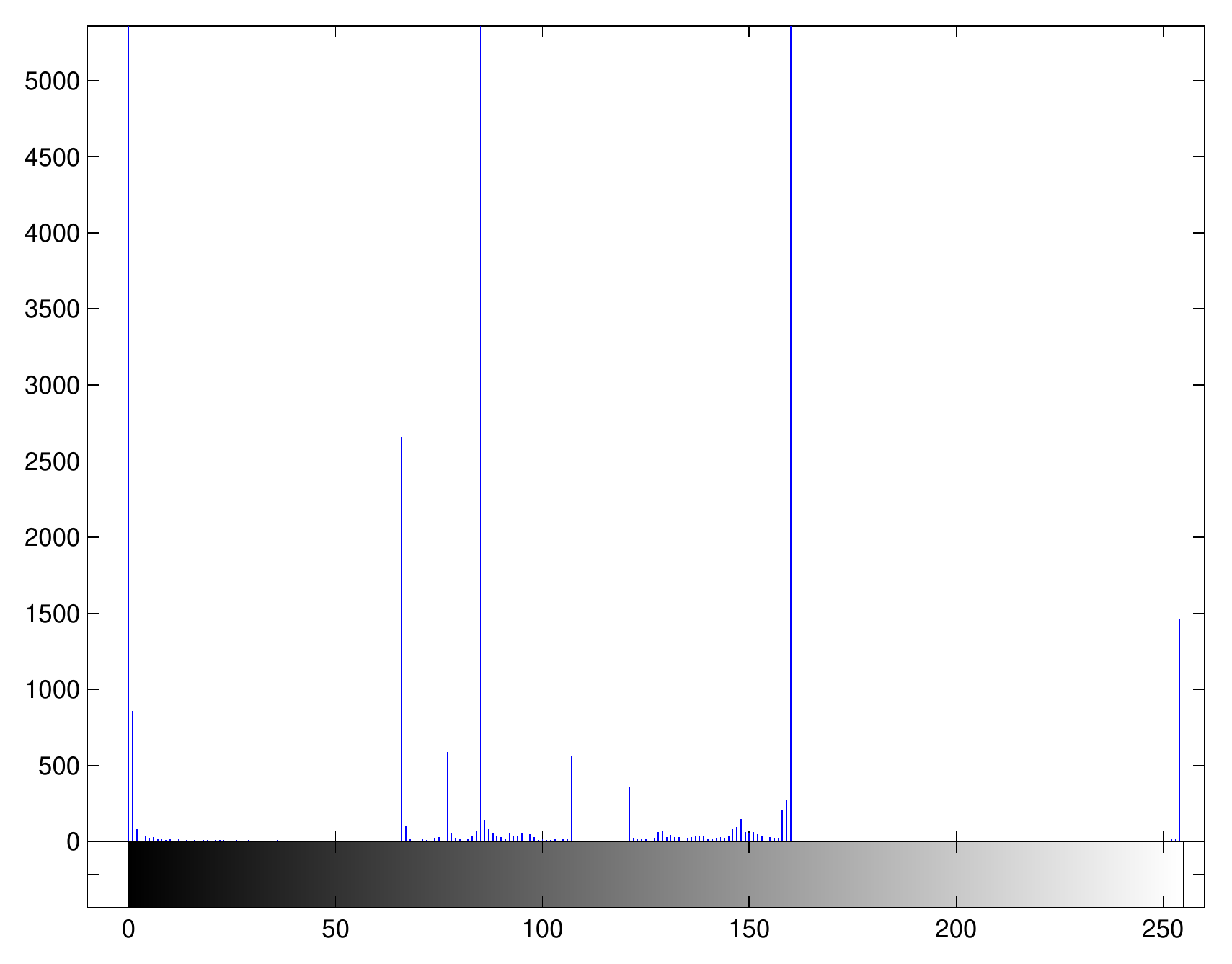} 
\end{minipage} 
\caption{The first 10 eigenvalues and the histograms of the first 3 eigenfunctions in Fig.\ref{fig:synth1}.
The $x$ and $y$ axes of the histograms are the gray value and the number of pixels having the same
gray value.} 
\label{fig:synth7hist} 
\end{figure}

The histograms of the first 3 eigenfunctions together with the plot
of the first 10 eigenvalues are shown in Fig.~\ref{fig:synth7hist}.
It is clear that the first 3 eigenfunctions are almost piecewise constant.
In fact, the fourth, fifth, and sixth are also almost piece constant whereas the seventh
is clearly not. (Their histograms are not shown here to save space but this can be seen
in Fig.~\ref{fig:synth1}.)

\begin{figure}
\begin{centering}
\includegraphics[width=10cm]{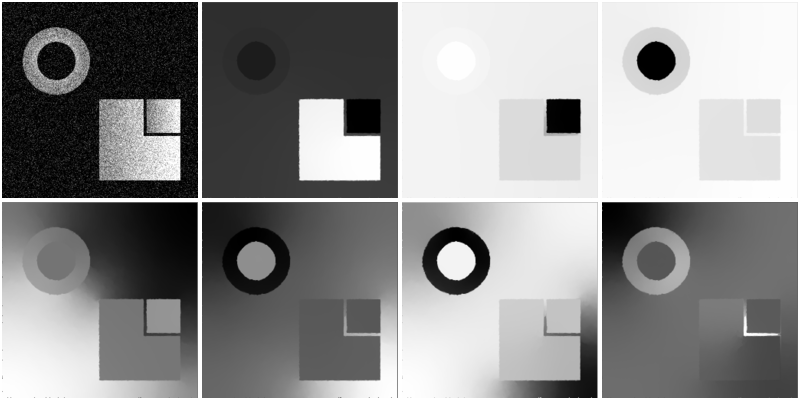}
\par\end{centering}
\caption{A noisy synthetic image and the first 7 eigenfunctions, left to right,
top to bottom. \label{fig:nzsynth1x}}
\end{figure}

Fig.\ref{fig:nzsynth1x} shows the results obtained an image with a mild level of noise.
The computation is done with the same condition as for Fig.~\ref{fig:synth1} except that
the input image is different. We can see that the first few eigenfunctions are also piecewise
constant and thus the phenomenon is relatively robust to noise.

\subsubsection{Eigenvalue problem (\ref{eq:HW-eigen}) versus Laplace-Beltrami operator}

Eigenfunctions of the Laplace-Beltrami operator (on surfaces) have been studied for image segmentation
\cite{Shah-00,Sochen-Kimmel-Malladi-98} and shape analysis \cite{Reuter-09,Reuter-06}.
Thus, it is natural to compare the performance of
the Laplace-Beltrami operator and that of the eigenvalue problem (\ref{eq:HW-eigen}).
For this purpose, we choose a surface such that the Laplace-Beltrami operator has the same
diffusion matrix as that defined in (\ref{D-2}), (\ref{D-4}), and (\ref{D-6}) and takes the form as
\begin{equation}
- \nabla \cdot \left(\frac{1}{\sqrt{1+|\nabla u|^{2}}}\begin{bmatrix}1+|\partial_y u_{0}|^{2} & -|\partial_x u_{0}\partial_y u_{0}|\\
-|\partial_x u_{0}\partial_y u_{0}| & 1+|\partial_x u_{0}|^{2}
\end{bmatrix}\nabla u\right) = \lambda \sqrt{1+|\nabla u_0|^{2}}\; u .
\label{LB-1}
\end{equation}
The main difference between this eigenvalue problem with (\ref{eq:HW-eigen}) is that there is a weight function on the right-hand side of (\ref{LB-1}), and 
in our model the parameter $\alpha$ in (\ref{D-6}) is typically greater than 1.

\begin{figure}
\begin{centering}
\includegraphics[width=10cm]{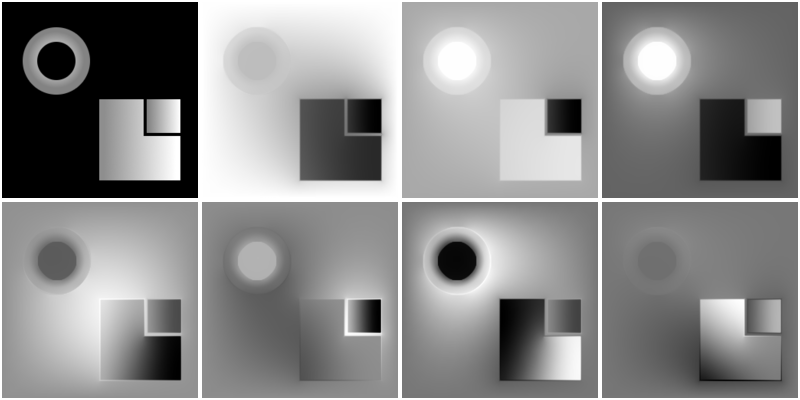}
\par\end{centering}
\caption{The clean input image of Fig.~\ref{fig:nzsynth1x} and the first 7
eigenfunctions of its associated Laplace-Beltrami operator.}
\label{fig:LB-eigenfunctions}
\end{figure}

The eigenfunctions of the Laplace-Beltrami operator obtained with a clean input image of Fig.~\ref{fig:nzsynth1x}
are shown in Fig.~\ref{fig:LB-eigenfunctions}. From these figures one can see that the eigenfunctions of
the Laplace-Beltrami operator are far less close to being piecewise constant, and thus, less suitable for
image segmentation.

\subsubsection{Open or closed edges}

\begin{figure}
\begin{centering}
\includegraphics[width=10cm]{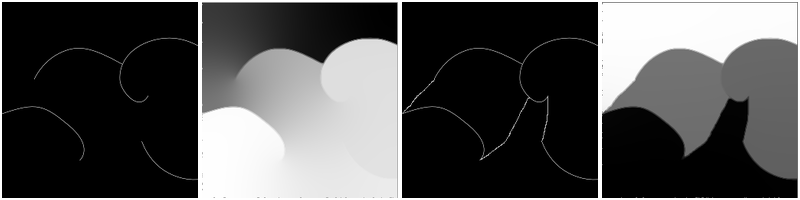}
\par\end{centering}
\caption{From left to right, input image with open arcs, the corresponding first eigenfunction,
input image with connected arcs, the corresponding first eigenfunction.}
\label{fig:openarc}
\end{figure}

We continue to study the piecewise constant property of eigenfunctions of (\ref{eq:HW-eigen}).
Interestingly, this property seems related to whether the edges of the input image form a closed curve.
We examine the two input images in Fig.~\ref{fig:openarc},  one containing
a few open arcs and the other having a closed curve that makes a jump in the gray level.
The first eigenfunction for the open-arc image changes gradually where 
that for the second image is close to being piecewise constant.

\subsubsection{Anisotropic mesh adaptation}

\begin{figure}
\begin{centering}
\includegraphics[width=10cm]{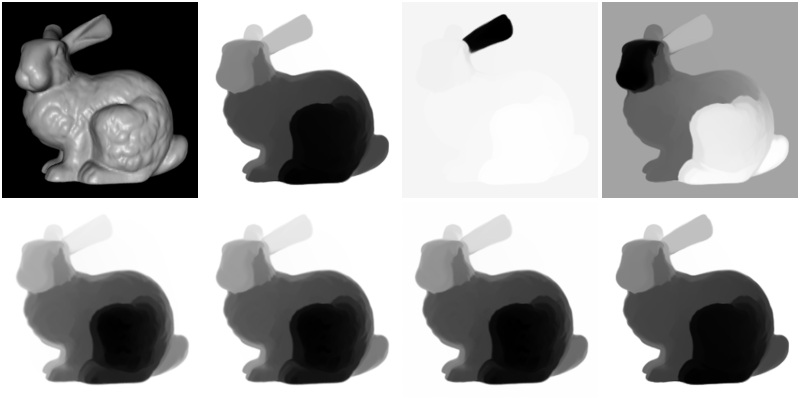}
\par\end{centering}
\caption{Top row: the image of the Stanford bunny and the first 3 eigenfunctions
computed with an anisotropic adaptive mesh with 45383 vertices; Bottom row: the
first eigenfunction on a uniform mesh with 93732, 276044, 550394 vertices and on
an adaptive mesh with 45383 vertices, respectively. All eigenfunctions are computed
with the same diffusion matrix defined in (\ref{D-2}), (\ref{D-4}), and (\ref{D-6}) ($\alpha = 1.5$)
and subject to the homogeneous Dirichlet boundary condition.}
\label{fig:bunny41}
\end{figure}

For the purpose of image segmentation, we would like the eigenfunctions to be as close to
being piecewise constant as possible. This would mean that they change abruptly in narrow regions
between objects. As a consequence, their numerical approximation can be difficult, and 
(anisotropic) mesh adaptation is then necessary in lieu of accuracy and efficiency.
The reader is referred to \cite{Huang-Wang-13} for the detailed studies of convergence and advantages
of using anisotropic mesh adaptation in finite element approximation of anisotropic eigenvalue problems
with anisotropic diffusion operators. Here, we demonstrate the advantage of using an anisotropic adaptive
mesh over a uniform one for the eigenvalue problem (\ref{eq:HW-eigen}) with the diffusion matrix defined
in (\ref{D-2}), (\ref{D-4}), and (\ref{D-6}) and subject to the homogeneous Dirichlet boundary condition.
The input image is taken as the Stanford bunny; see Fig.~\ref{fig:bunny41}. The figure also shows
the eigenfunctions obtained on an adaptive mesh and uniform meshes of several sizes.
It can be seen that the eigenfunctions obtained with the adaptive mesh have very sharp
boundaries, which are comparable to those obtained with a uniform mesh of more than
ten times of vertices.

\subsubsection{Anisotropic and less anisotropic diffusion}

\begin{figure}
\begin{centering}
\includegraphics[width=10cm]{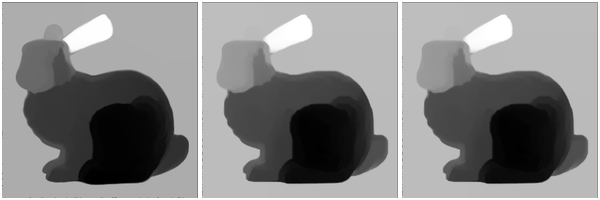}
\par\end{centering}
\caption{The first eigenfunction of (\ref{eq:HW-eigen}) with the Stanford bunny and
with the diffusion matrix (\ref{D-2}) and a less anisotropic one (\ref{D-7}).
From left to right,  (\ref{D-2}) on an adaptive mesh of 45562 vertices,
(\ref{D-7}) on a uniform mesh of size $256\times256$, and
(\ref{D-7}) on an adaptive mesh of 46974 vertices.}
\label{fig:ani-vs-iso}
\end{figure}

Next, we compare the performance of the diffusion matrix (\ref{D-2}) (with 
(\ref{D-4}), (\ref{D-6}), and $\alpha = 1.5$) and that of a less anisotropic diffusion matrix
(cf. (\ref{eq:pm-linear-diffusion}), with (\ref{D-6}) and $\alpha = 1.5$)
\begin{equation}
\mathbb{D} = \begin{bmatrix}g(|\partial_{x}u_0 |) & 0\\ 0 & g(|\partial_{y} u_0|) \end{bmatrix} .
\label{D-7}
\end{equation}
The eigenfunctions of (\ref{eq:HW-eigen}) with those diffusion matrices
with the Stanford bunny as the input image are shown in 
Fig.~\ref{fig:ani-vs-iso}. For (\ref{D-7}), we compute the
eigenfunction on both a uniform mesh of size $256\times256$ and
an adaptive mesh of 46974 vertices. The computation with (\ref{D-2})
is done with an adaptive mesh of 45562 vertices.
The most perceptible difference in the results is that the right ear of the bunny
(not as bright as other parts) almost disappears in the first eigenfunction with
the less anisotropic diffusion matrix. This can be recovered if the conductance
is changed from $\alpha = 1.5$ to $\alpha = 1.0$, but in this case,
the eigenfunction becomes farther from being piecewise-constant.
The image associated with the first eigenfunction for (\ref{D-2}) seems sharper
than that with (\ref{D-7}). 

\subsubsection{Effects of the conductance $g$}
\label{SEC:4.1.6}

\begin{figure} 
\begin{minipage}{1.\linewidth} 
\centering 
\includegraphics[scale=.4]{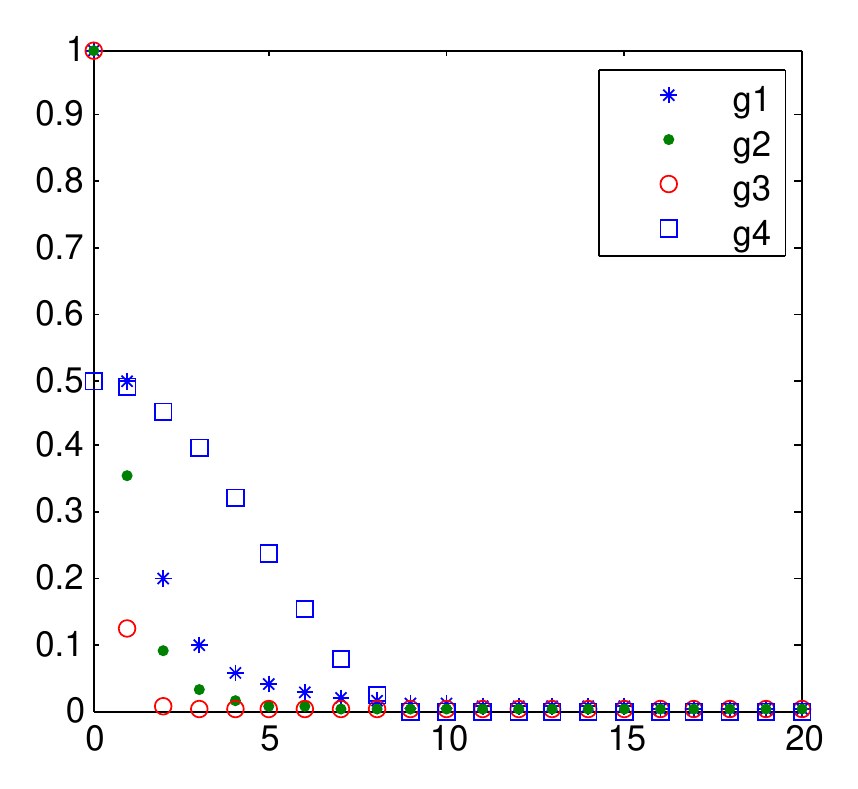}
\end{minipage} \par\medskip 

\begin{minipage}{.2\linewidth} 
\centering 
\includegraphics[scale=.3]{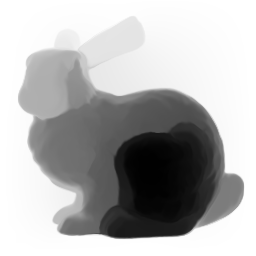}
\end{minipage}%
\hfill
\begin{minipage}{.2\linewidth} 
\centering 
\includegraphics[scale=.3]{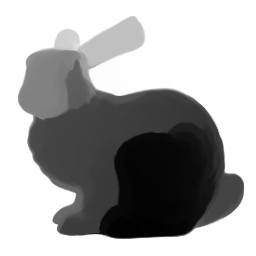} 
\end{minipage} 
\hfill
\begin{minipage}{.2\linewidth} 
\centering 
\includegraphics[scale=.3]{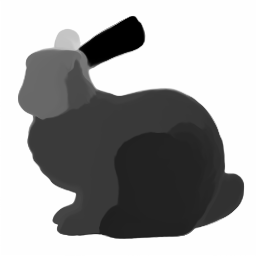} 
\end{minipage} 
\hfill
\begin{minipage}{.2\linewidth} 
\centering 
\includegraphics[scale=.3]{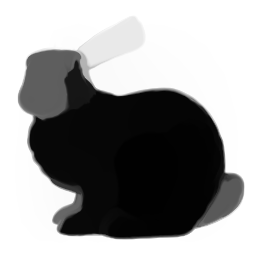} 
\end{minipage} \par\medskip 

\begin{minipage}{.2\linewidth} 
\centering 
\includegraphics[scale=.35]{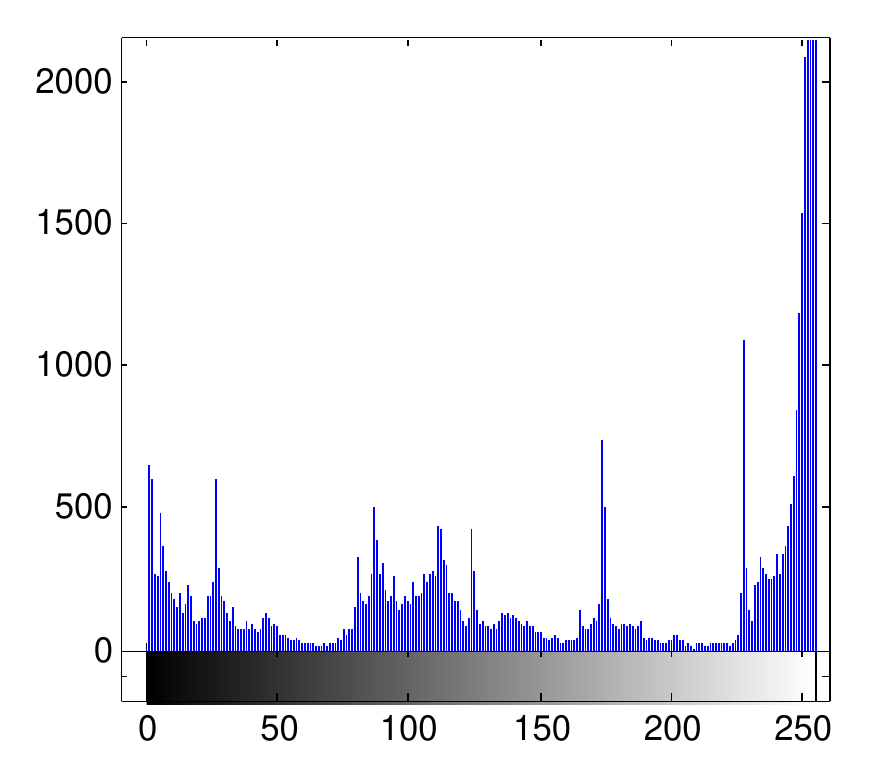}
\end{minipage}%
\hfill
\begin{minipage}{.2\linewidth} 
\centering 
\includegraphics[scale=.35]{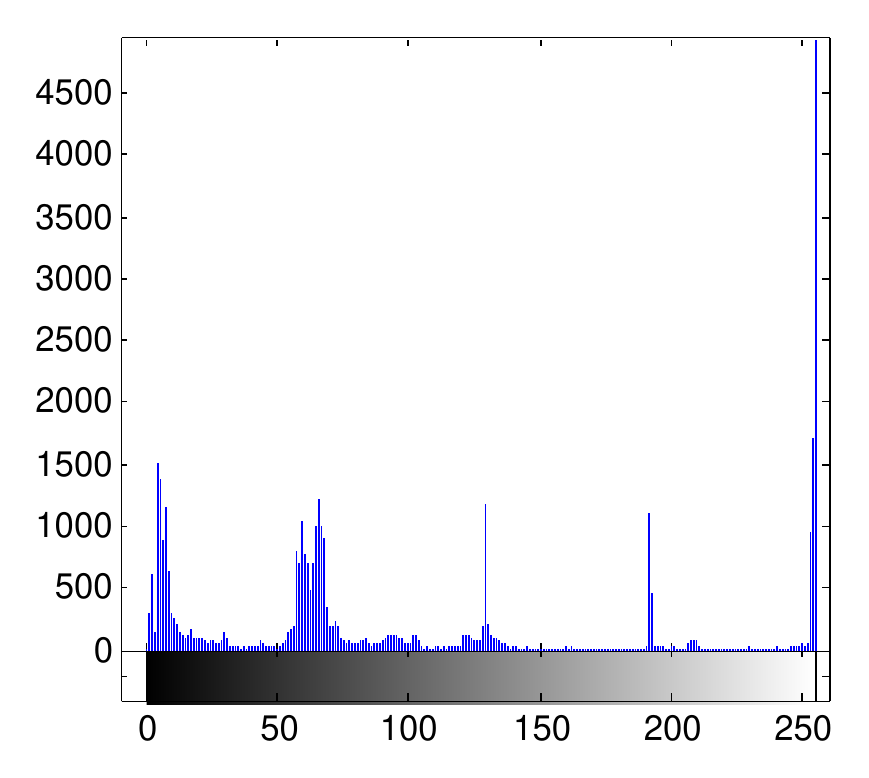} 
\end{minipage} 
\hfill
\begin{minipage}{.2\linewidth} 
\centering 
\includegraphics[scale=.35]{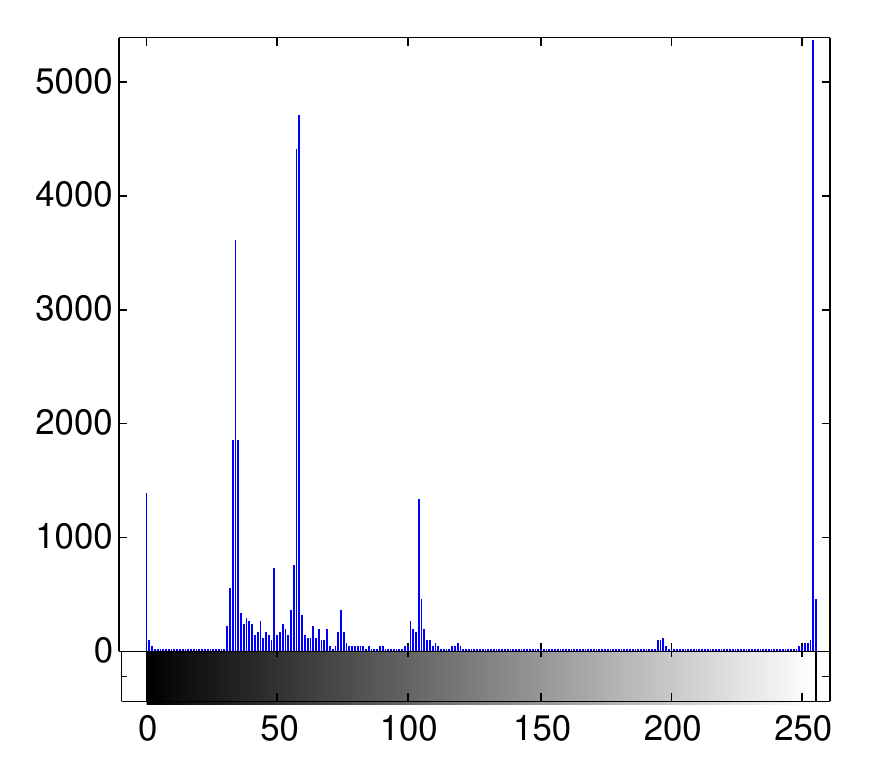} 
\end{minipage} 
\hfill
\begin{minipage}{.2\linewidth} 
\centering 
\includegraphics[scale=.35]{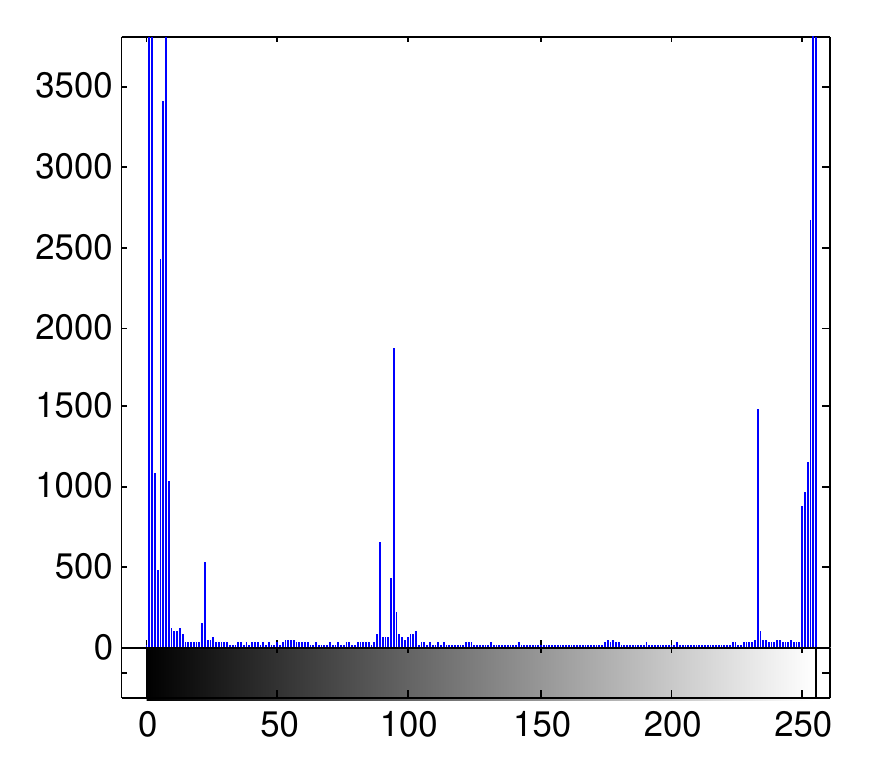} 
\end{minipage} 
\caption{Top row: the graphs of $g_1, g_2, g_3, g_4$. Middle row: the first eigenfunctions on the bunny image for $g_1, g_2, g_3, g_4$, respectively. Bottom row: the histograms of the corresponding first eigenfunctions.} 
\label{fig:g-choice} 
\end{figure}

We now examine the effects of the conductance and consider four cases:
$g_1$ ((\ref{D-6}) with $\alpha = 1.0$), $g_2$ ((\ref{D-6}) with $\alpha = 1.5$),
$g_3$ ((\ref{D-6}) with $\alpha = 3.0$), and
\[
g_{4}(x)=\begin{cases}
(1-(x/\sigma)^{2})^{2}/2, & \text{ for }|x|\le\sigma\\
\epsilon, & \text{ for }|x|>\sigma
\end{cases}
\]
where $\sigma$ and $\epsilon$ are positive parameters.
The last function is called Tukey's biweight function and considered
in \cite{Black-Sapiro-01} as a more robust choice of the edge-stopping
function in the Perona-Malik diffusion. We show the results with
(\ref{D-2}) on the Stanford bunny in Fig.~\ref{fig:g-choice}.
We take $\sigma=9$ and $\epsilon=10^{-6}$ for Tukey's biweight function. Increasing
the power $\alpha$ in $g(x)$ defined in (\ref{D-6}) will make
eigenfunctions steeper in the regions between different objects
and thus, closer to being piecewise constant. Tukey's biweight function
gives a sharp result but the body and legs are indistinguishable.

\subsection{Applications in Edge Detection and Image Segmentation}

Eigenfunctions can serve as a low level image feature extraction device
to facilitate image segmentation or object edge detection. Generally speaking,
eigenfunctions associated with small eigenvalues contain ``global'' segmentation features
of an image while eigenfunctions associated with larger eigenvalues carry more information
on the detail. Once the eigenfunctions are obtained, one can use
numerous well developed edge detection or data clustering techniques
to extract edge or segmentation information. We point out that spectral
clustering methods also follow this paradigm. In this section,
we focus on the feature extraction step and employ only simple, well known
techniques such as thresholding by hand, $k$-means
clustering, or Canny edge detector in the partitioning step. More
sophisticated schemes can be easily integrated to automatically detect
edges or get the segmentations.

We point out that boundary conditions have an interesting effect on
the eigenfunctions. A homogeneous Dirichlet boundary condition
forces the eigenfunctions to be zero on the boundary and may wipe out
some structures there (and therefore, emphasize objects inside the domain).
It essentially plays the role of defining ``seeds'' that indicates 
background pixels on the image border. The idea of using user-defined seeds or 
intervene cues has been widely used in graph based image segmentation methods 
\cite{Grady-01}, \cite{Rother-Blake-04}, \cite{Shi-Yu-04}, \cite{Malik-Martin-04}. 
The PDE eigenvalue problem (\ref{eq:HW-eigen}) can also be solved with more sophisticated 
boundary conditions that are defined either on the image border or inside the image.
On the other hand, a homogeneous Neumann
boundary condition tends to keep those structures. 
Since mostly we are interested in objects inside the domain, we consider here
a homogeneous Dirichlet boundary condition. The diffusion matrix
defined in (\ref{D-2}), (\ref{D-4}), and (\ref{D-6}) ($\alpha = 1.5$) is used.

\begin{figure}
\begin{centering}
\includegraphics[width=10cm]{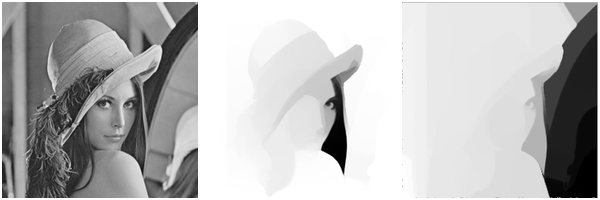}
\par\end{centering}
\caption{From left to right, Lenna, first eigenfunctions obtained with Dirichlet and Neumann
boundary conditions, respectively.}
\label{fig:Lenna_ef1}
\end{figure}

\begin{figure}
\centering{}\includegraphics[width=12cm]{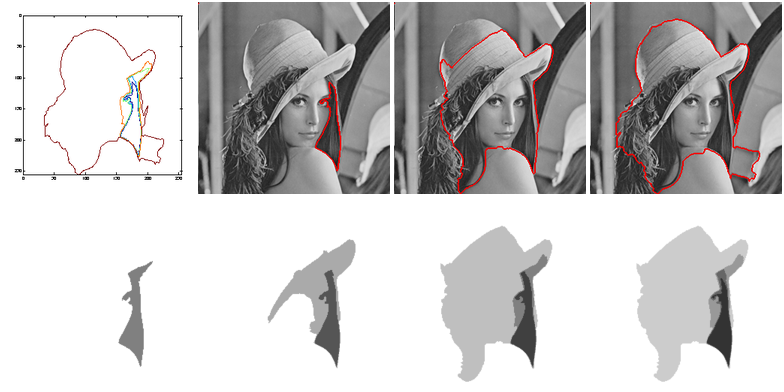}
\caption{Top row: contour drawing by MATLAB (with no level parameters specified),
level line 50, 240, 249; bottom row: segmentation with $k$-means, $k=2,3,4,5$.}
\label{fig:lenna-contour}
\end{figure}

In Fig.~\ref{fig:Lenna_ef1}, we show the first eigenfunctions obtained with Dirichlet and Neumann
boundary conditions with Lenna as the input image.
For the edge detection for Lenna, it is natural to extract the ``big picture'' from the first eigenfunction
and get the edge information from it. We show the edges obtained by thresholding a few level lines
in the top row of Fig.~\ref{fig:lenna-contour}. Since any level line
with value $s$ is the boundary of the level set $L_{s}=\{(x,y):I(x,y)\ge s\}$
of an image $I$, and $L_{s}$ is non-increasing with respect to $s$,
the level line is ``shrinking'' from the boundary of a wider shape
to empty as $s$ increases from 0 to 255. Some intermediate steps
give salient boundaries of the interior figure. However, to make the ``shrinking''
automatically stop at the correct edge, other clues potentially from
mid or high level knowledge in addition to the low level brightness
info should be integrated in the edge detection step. We also use
the MATLAB function {\tt imcontour} to get major contours, and apply $k$-means
clustering to the eigenfunctions with $k=2,3,4,5$, shown in the second
row of Fig.~\ref{fig:lenna-contour}.

\begin{figure}
\begin{centering}
\includegraphics[width=14cm]{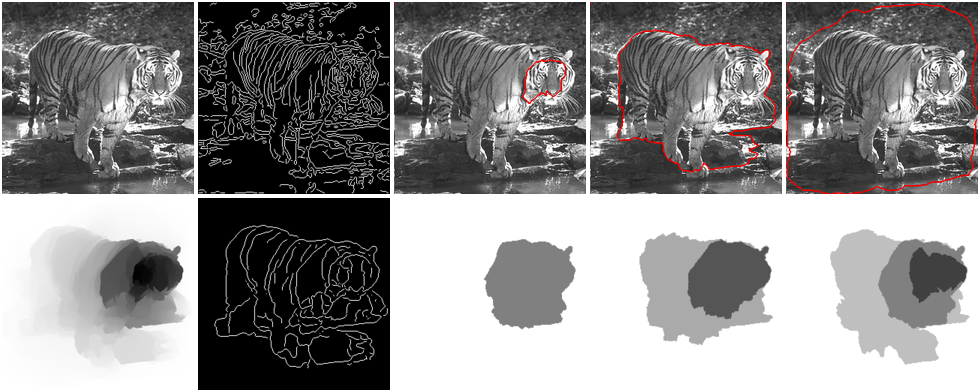}
\par\end{centering}
\caption{Top row: the input image, the edges of the input image with the Canny detector,
Level lines with value 50, 240, 254. bottom row: the first eigenfunction,
the edges of the first eigenfunction with the Canny detector, $k$-means
clustering with $k=2,3,4$.}
\label{fig:tiger-gallery}
\end{figure}

We next compute for an image with more textures from \cite{MartinFTM01} (Fig.~\ref{fig:tiger-gallery}).
This is a more difficult image for
segmentation or edge detection due to many open boundary arcs
and ill-defined boundaries. We display the the first eigenfunction
and the $k$-means clustering results in Fig.~\ref{fig:tiger-gallery}.
The $k$-means clustering does not capture the object as well as in
the previous example. Better separation of the object and the background
can be obtained if additional information is integrated into the clustering strategy.
For instance, the edges detected by the Canny detector (which uses
the gradient magnitude of the image) on the
first eigenfunction clearly give the location of the tiger. Thus,
the use of the gradient map of the first eigenfunction
in the clustering process yields more accurate object boundaries.
For comparison, we also show the edges detected from the input image
with the Canny detector.

\begin{figure} 
\begin{minipage}{.2\linewidth} 
\centering 
\includegraphics[scale=.3]{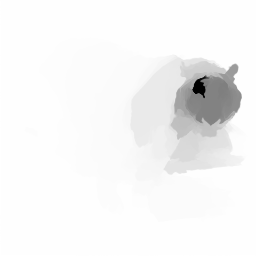}
\end{minipage}%
\hfill
\begin{minipage}{.2\linewidth} 
\centering 
\includegraphics[scale=.3]{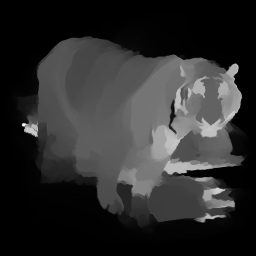} 
\end{minipage} 
\hfill
\begin{minipage}{.2\linewidth} 
\centering 
\includegraphics[scale=.35]{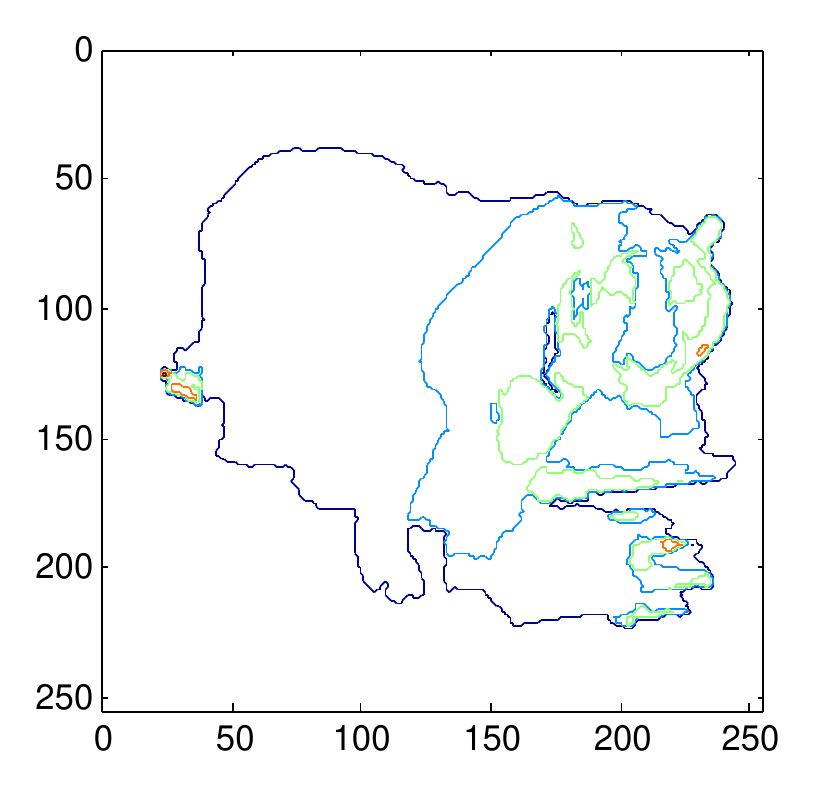} 
\end{minipage} 
\hfill
\begin{minipage}{.2\linewidth} 
\centering 
\includegraphics[scale=.35]{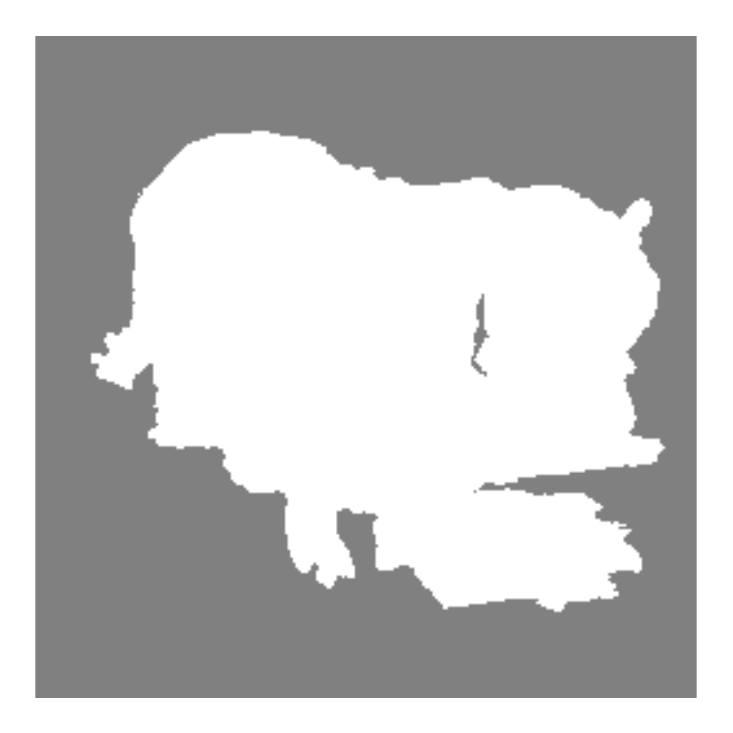} 
\end{minipage} \par\medskip 
\caption{From left to right, the first eigenfunction with $\alpha = 2$ in (\ref{D-6}),
the projection of the tiger image into the space spanned by the first 64 eigenfunctions, the contour of the projection
image, the $k$-means clustering with $k=2$ of the projection image.} 
\label{fig:tigerfin} 
\end{figure}

Another way to extract ``simple'' features is to change the conductance
$g$ (e.g., by increasing $\alpha$ in (\ref{D-6})) to make
the eigenfunctions closer to being piecewise constant. This makes eigenfunctions more clustered
but wipes out some detail of the image too. To avoid this difficulty,
we can employ a number of eigenfunctions and use the projection of
the input image into the space spanned by the eigenfunctions to construct
a composite image. A much better result obtained in this way with 64 eigenfunctions
is shown in Fig.~\ref{fig:tigerfin}.

\begin{figure}
\begin{centering}
\includegraphics[width=14cm]{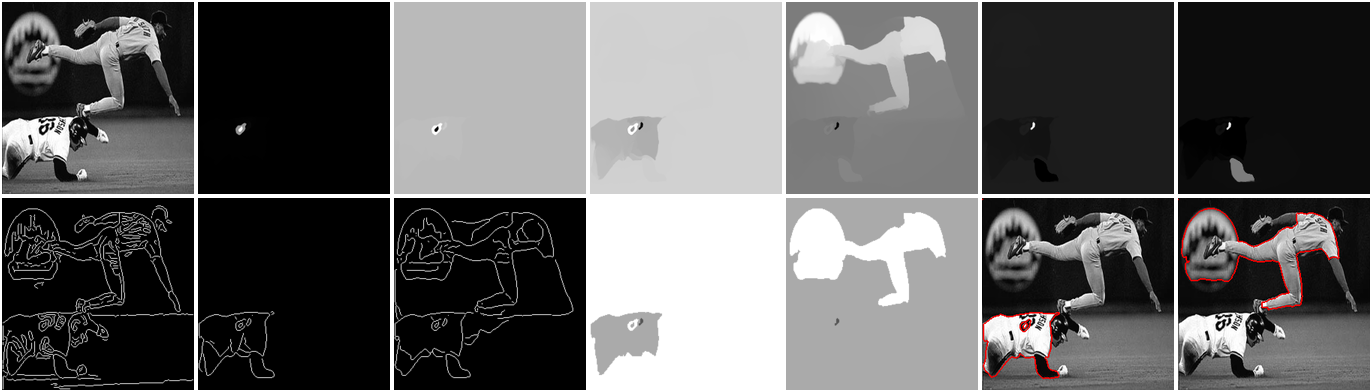}
\par\end{centering}
\caption{Top row: from left to right, the input image and the first 6 eigenfunctions.
Bottom row: from left to right, the edges on the input image (Canny),
the edges on the 3rd and 4th eigenfunctions (Canny), the $k$-means
clustering results with $k=3$ for the 3rd and the 4th eigenfunctions,
respectively; Level line of value 205 of the 3rd eigenfunction, level
line of value 150 of the 4th eigenfunction.}
\label{fig:sports-gallery}
\end{figure}

\begin{figure}
\begin{centering}
\includegraphics[scale=0.6]{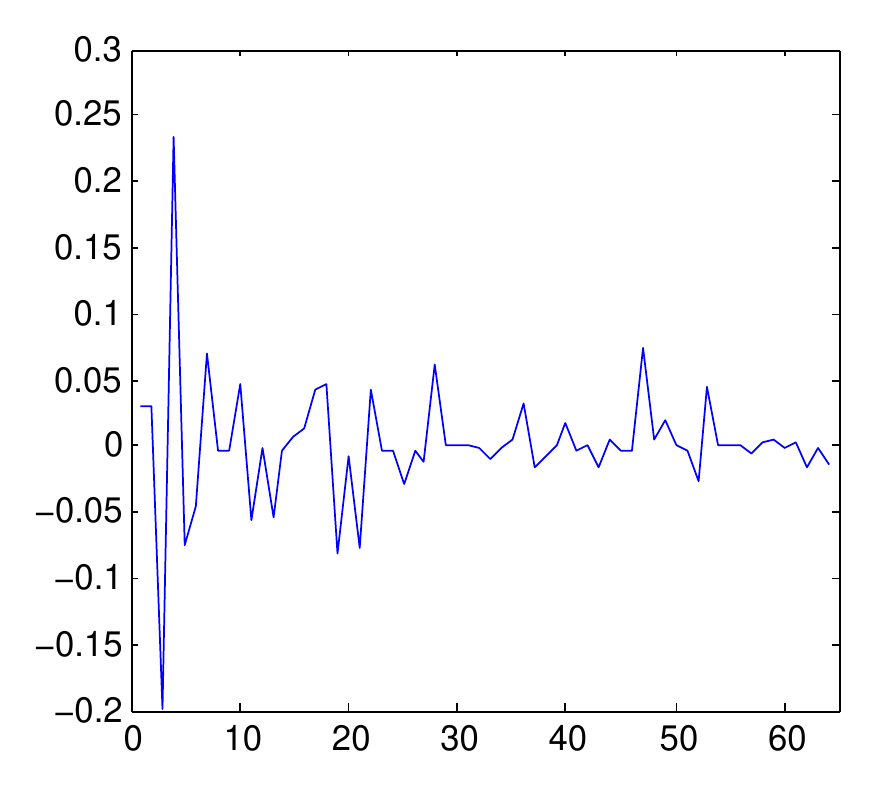}
\par\end{centering}
\caption{The coefficients of the input image projected onto the first 64 eigenfunctions in
Fig.~\ref{fig:sports-gallery}.}
\label{fig:sports-components}
\end{figure}

It should be pointed out that not always the first few eigenfunctions cary most useful information
of the input image. Indeed, Fig.~\ref{fig:sports-gallery} shows that
the first eigenfunction carries very little information. Since the
eigenfunctions form an orthogonal set in $L^{2}$, we can project
the input image onto the computed eigenfunctions.
The coefficients are shown in Fig.~\ref{fig:sports-components}.
We can see that the coefficients for the first two eigenfunctions are very small
compared with those for the several following eigenfunctions.
It is reasonable to use the eigenfunctions with the greatest magnitudes of
the coefficients. These major eigenfunctions will provide most useful information;
see Fig.~\ref{fig:sports-gallery}.

\section{The piecewise constant property of eigenfunctions}
\label{SEC:piecewise}

As we have seen in Section~\ref{SEC:numerics}, 
the eigenfunctions of problem (\ref{eq:HW-eigen})
are localized in sub-regions of the input image and the first few of them are close
to being piecewise constant for most input images except for two types of images.
The first type of images is those containing
regions of which part of their boundaries is not clearly defined (such
as open arcs that are common in natural images). In this case, the first eigenfunction
is no longer piecewise-constant although the function values can still be well clustered.
The other type is input images for which the gray level changes gradually and its gradient
is bounded (i.e., the image contrast is mild).
In this case, the diffusion operator simply behaves like the Laplace operator and has
smooth eigenfunctions. For other types of images, the gray level has an abrupt change
across the edges of objects, which causes the conductance $g(|\nabla u_0|)$
to become nearly zero on the boundaries between the objects. As a consequence,
the first few eigenfunctions are close to being constant within each object. This property
forms the basis for the use of the eigenvalue problem (\ref{eq:HW-eigen}) (and its eigenfunctions)
in image segmentation and edge detection.
In this section, we attempt to explain this property from the physical, mathematical, and graph spectral points of view.
We hope that the analysis, although not rigorous, provides some insight of the phenomenon.

From the physical point of view, when the conductance $g(|\nabla u_0|)$
becomes nearly zero across the boundaries between the objects, the diffusion flux will be nearly zero
and each object can be viewed as a separated region from other objects. As a consequence,
the eigenvalue problem can be viewed as a problem defined on multiple separated subdomains, subject to
homogeneous Neumann boundary conditions (a.k.a. insulated boundary conditions)
on the boundary of the whole image and the internal 
boundaries between the objects. Then, it is easy to see that the eigenfunctions corresponding to the eigenvalue 0
include constant and piecewise constant (taking a different constant value on each object) functions.
This may explain why piecewise constant eigenfunctions have been observed for most input images.
On the other hand, for images with mild contrast or open arc object edges, the portion of the domain associated
any object is no longer totally separated from other objects and thus the eigenvalue problem may not have piecewise constant
eigenfunctions.

\begin{figure}
\begin{centering}
\includegraphics[scale=0.4]{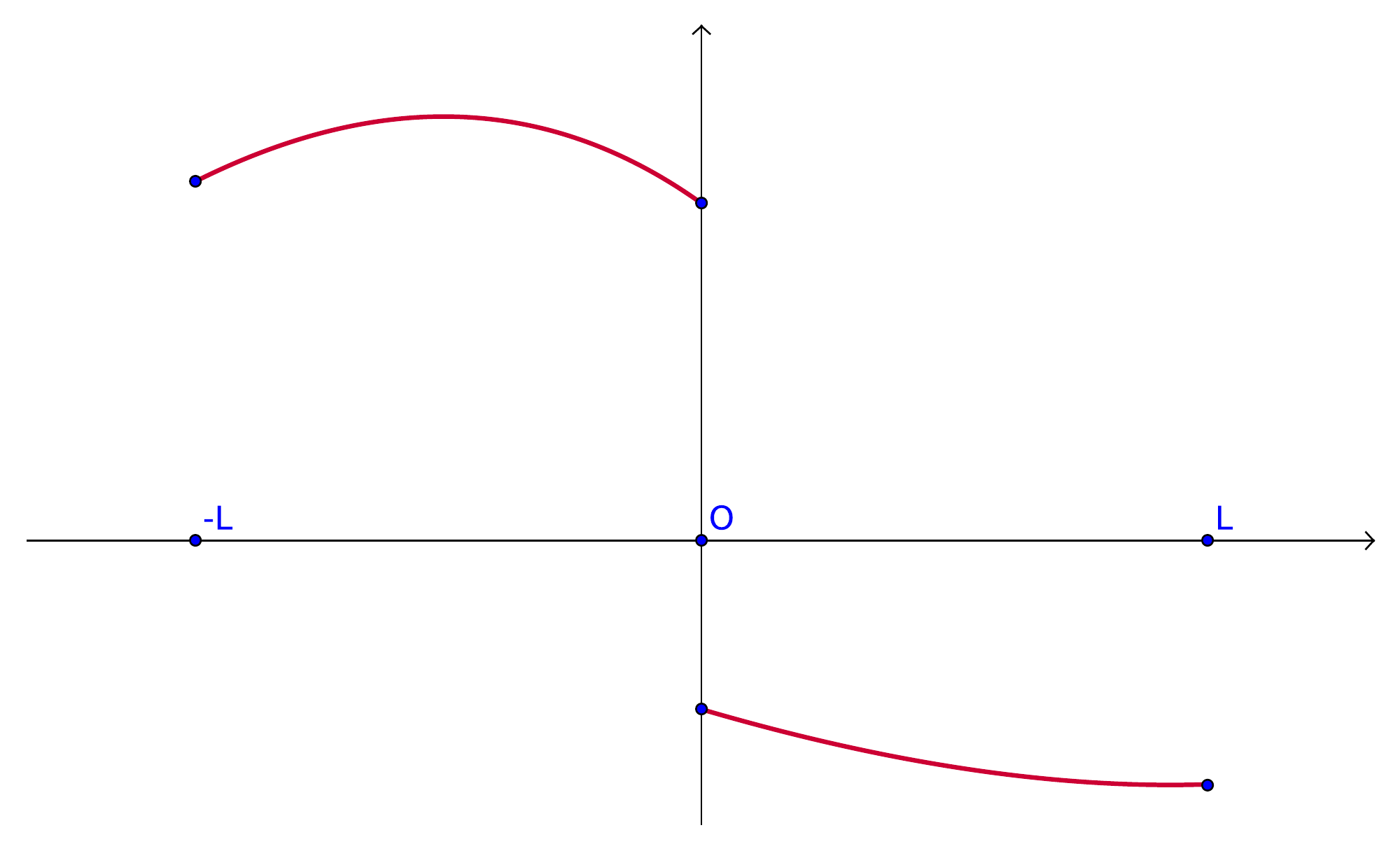}
\par\end{centering}
\caption{A piecewise smooth function representing an input image with two objects.}
\label{fig:a-fun}
\end{figure}

Mathematically, we consider a 1D example with an input image with two objects. The gray level of the image is sketched
in Fig.~\ref{fig:a-fun}. The edge is located at the origin, and the
segmentation of the image is a 2-partition of $[-L,0]$ and
$[0,L]$. The 1D version of the eigenvalue problem (\ref{eq:HW-eigen}) reads as
\begin{equation}
- \frac{d}{d x}\left(g(|u_0'(x) |)\frac{d u}{d x} \right) =\lambda u ,
\label{eq:example1}
\end{equation}
subject to the Neumann boundary conditions $u'(-L)=u'(L)=0$. We take the conductance
function as in (\ref{D-6}) with $\alpha = 2$. Although $u_0$ is not differentiable,
we could imagine that $u_0$ were replaced by a smoothed function which has a very large
or an infinite derivative at the origin. Then, (\ref{eq:example1}) is degenerate since $g(|u_0'(x)|$
vanishes at $x=0$. As a consequence, its eigenfunctions can be non-smooth. Generally speaking,
studying the eigenvalue problem of a degenerate elliptic operator is a difficult task,
and this is also beyond the scope of the current work. Instead of performing a rigorous analysis, we
consider a simple approximation to $g(|u_0'(x)|)$,
\[
g_{\epsilon}(x)=\begin{cases}
g(|u_0'(x)|), & \text{ for }-L\le x<-\epsilon{\rm \ or\ }\epsilon<x\le L\\
0, & \text{ for }-\epsilon \le x\le \epsilon
\end{cases}
\]
where $\epsilon$ is a small positive number.
The corresponding approximate eigenvalue problem is
\begin{equation}
- \frac{d}{d x}\left(g_\epsilon(|u_0'(x) |)\frac{d u}{d x} \right) =\lambda u .
\label{eq:exampe1-approx}
\end{equation}
The variational formulation of this eigenvalue problem is given by
\begin{equation}
\min_{u\in H^1(-L,L)}\int_{-L}^{L}g_{\epsilon}(|u_0'(x)|)(u')^{2}, \quad \text{ subject to } \int_{-L}^{L}u^{2}=1 .
\label{eq:ex1-variation-prob}
\end{equation}
Once again, the eigenvalue problem (\ref{eq:exampe1-approx}) and the variational problem (\ref{eq:ex1-variation-prob})
are degenerate and should be allowed to admit non-smooth solutions.

The first eigenvalue of (\ref{eq:exampe1-approx}) is 0, and a trivial eigenfunction associated with this eigenvalue
is a constant function. To get other eigenfunctions associated with 0, we consider
functions that are orthogonal to constant eigenfunctions, i.e., we append to
the optimization problem (\ref{eq:ex1-variation-prob}) with the constraint
\[
\int_{-L}^{L}u=0.
\label{eq:ex1-constraints}
\]
It can be verified that an eigenfunction is
\begin{equation}
u_{\epsilon}(x)=\begin{cases}
-c, & \text{ for }x\in [-L,-\epsilon)\\
\frac{c x}{\epsilon}, & \text{ for } x\in [-\epsilon, \epsilon]\\
c, & \text{ for } x\in (\epsilon, L]
\end{cases}
\label{1D-solution}
\end{equation}
where $c = (2 (L-2\epsilon/3))^{-1/2}$. This function is piecewise constant for most part of the domain
except the small region $[-\epsilon, \epsilon]$. Since the original problem (\ref{eq:example1}) can be viewed
to some extent as the limit of (\ref{eq:exampe1-approx}) as $\epsilon \to 0$, the above analysis may
explain why some of the eigenfunctions of (\ref{eq:example1}) behave like piecewise constant functions.

The piecewise constant property can also be understood in the context of the graph spectral
theory. We first state a result from \cite{Mohar-Alavi-91,Luxburg-2007}.

\begin{prop}[\cite{Mohar-Alavi-91,Luxburg-2007}]\; 
Assume that $G$ is a undirected graph with $k$ connected components and
the edge weights between those components are zero. If the nonnegative weights matrix $W$ 
and the diagonal matrix $D$ are defined as in Section~\ref{SEC:relation-1},
then the multiplicity of the eigenvalue 0 of the matrix $D-W$ equals the number of the connected
components in the graph. Moreover, the eigenspace of
eigenvalue 0 is spanned by the indicator vectors $1_{A_{1}},\cdots,1_{A_{k}}$
of those components ${A_{1}},\cdots, {A_{k}}$.
\label{prop-1}
\end{prop}

This proposition shows that the eigenvalue zero of the algebraic eigenvalue problem (\ref{eq:Malik-Shi-eigen}) 
could admit multiple eigenvectors as indicators of the components. As shown in Section~\ref{SEC:relation-1},
(\ref{eq:Malik-Shi-eigen}) can be derived from a finite difference discretization of the continuous
eigenvalue problem (\ref{eq:HW-eigen}) (with a proper choice of $\mathbb{D}$). Thus,
the indicators of the components can also be regarded as discrete approximations of some continuous
eigenfunctions. This implies that the latter must behave like piecewise constant functions.
Interestingly, Szlam and Bresson \cite{Szlam-Bresson-10} recently proved
that global binary minimizers exist for a graph based problem called Cheeger Cut where the minimum of
the cut is not necessarily zero.

In the continuous setting, a properly designed conductance $g(|\nabla u_0|)$
can act like cutting the input image $u_0$ into subregions along the boundary
of the subregions and forcing the eigenvalue problem to be solved
on each subregion. In a simplified case, we can have the following continuous
analogue of Proposition~\ref{prop-1}. The proof is straightforward.

\begin{prop}\;
Suppose $\Omega\subset \mathbb{R}^2 $ is a bounded Lipschitz domain, $u_0 \in SBV(\Omega)$
(the collect of special functions of bounded variation) and the discontinuity set $K$ of $u_0$ is a finite
union of $C^{1}$ closed curves, $g(\cdot)$ is a bounded positive continuous function. 
We define $g(|\nabla u_0|)=0$
for $(x,y)\in K$, and $g(|\nabla u_0|)$ takes its usual meaning for $(x,y)\in\Omega\backslash K$.
For any function $u\in SBV(\Omega)$, assuming $\Gamma$ is the discontinuity
set of $u$, we define the Rayleigh quotient on $u$ as
\[
R(u)=\frac{\widetilde{\int}_{\Omega}(\nabla u)^{T}g(|\nabla u_0|)\nabla u}{\int_{\Omega}u^{2}} ,
\]
where
\[
\widetilde{\int}_{\Omega}(\nabla u)^{T}g(|\nabla u_0|)\nabla u=\begin{cases}
\int_{\Omega\backslash K}g(\nabla u_0|)|\nabla u|^{2}, & \text{ for }\Gamma\subseteq K\\
\infty, & \text{ for } \Gamma\nsubseteq K .
\end{cases}
\]
Then, the minimum of $R(u)$ is zero and any piecewise constant function
in $SBV(\Omega)$ with discontinuity set in $K$ is a minimizer.
\label{prop-2}
\end{prop}

%
%

The eigenvalue problem related to the above variational problem can be formally written as

\begin{equation}
-\nabla \cdot \left(g(|\nabla u_0|)\nabla u\right)=\lambda u.
\label{eq:general-eigen-prob}
\end{equation}
The equation should be properly defined for all $u_0$, $u$ possibly in the space
of $BV$. This is a degenerate elliptic problem which could admit discontinuous
solutions, and it seems to be far from being fully understood.
In the following proposition, we suggest a definition of weak solution
in a simplified case. The property indicates that problem (\ref{eq:general-eigen-prob})
is quite different from a classical elliptic eigenvalue problem if
it has a solution in $BV$ that takes a non-zero constant value on an
open set. The proof is not difficult and thus omitted.
\begin{prop}\;
Suppose $\Omega$ is a bounded Lipschitz domain in $R^{2}$, $u_0\in SBV(\Omega)$
and the discontinuity set $K$ of $u_0$ is a finite union of $C^{1}$
closed curves, $g(\cdot)$ is a bounded positive continuous function.
We define $g(|\nabla u_0|)=0$ for $(x,y)\in K$, and $g(|\nabla u_0|)$ takes its usual 
meaning for $(x,y)\in\Omega\backslash K$.
We define $u\in SBV(\Omega)$ to be a weak eigenfunction of (\ref{eq:general-eigen-prob})
satisfying a homogeneous Dirichlet boundary condition if
\[
\int_{\Omega}(\nabla u)^{T}g(|\nabla u_0|)\nabla\phi=\int_{\Omega}\lambda u\phi,
\quad \forall \phi \in C_{0}^{1}(\Omega)
\]
where, assuming that $\Gamma$ is the discontinuity set of $u$, the integral
on the left side is defined by
\[
\int_{\Omega}(\nabla u)^{T}g(|\nabla u_0|)\phi=\begin{cases}
\int_{\Omega\backslash K}(\nabla u)^{T}g(\nabla u_0|)\phi, & \Gamma\subseteq K\\
\infty, & {\rm \Gamma\nsubseteq K.}
\end{cases}
\]
If a weak eigenfunction $u\in SBV(\Omega)$ exists and takes a non-zero constant value
on a ball $B_{\epsilon}(x_{0}, y_0)\subset\Omega$, then the corresponding eigenvalue $\lambda$ is zero.
\label{prop-3}
\end{prop}

If (\ref{eq:general-eigen-prob}) indeed admits non-zero
piecewise-constant eigenfunctions, one can see an interesting connection
between (\ref{eq:general-eigen-prob}) (for simplicity we assume
a homogeneous Dirichlet boundary condition is used) and Grady's Random Walk image
segmentation model \cite{Grady-01} where multiple combinatorial Dirichlet problems are solved
for a $k$-region segmentation with predefined seeds indicating segmentation labels.
Using a similar argument in Section 2.2, one can show that the numerical
implementation of the method is equivalent to solving a set of Laplace problems
which are subject to a Neumann boundary condition on the image border and Dirichlet boundary conditions
on the seeds and are discretized on a uniform mesh for potentials $u^{i}$,
$i=1,\cdots,k$. These boundary problems read as
\begin{eqnarray}
\nabla \cdot \left(\begin{bmatrix}g(|\partial_{x}u_{0}|) & 0\\
0 & g(|\partial_{y}u_{0}|)
\end{bmatrix}\nabla u^{i}\right) & = & 0,\ {\rm in}\ \Omega\backslash S
\label{eq:Grady-model}\\
\frac{\partial u^{i}}{\partial n} & = & 0,\ {\rm on}\ \partial\Omega \nonumber \\
u^{i} & = & 1,\ {\rm in}\ S_{i} \nonumber \\
u^{i} & = & 0,\ {\rm in}\ S\backslash S_{i}
\nonumber 
\end{eqnarray}
where $S_{i}$ is the set of seeds for label $i$ and $S$ is the set of all seeds.
This problem with a proper choice of $g$ also gives
a solution that has well clustered function values, a phenomenon called
``histogram concentration'' in \cite{Buades-Chien-08}. 

Note that when $\lambda=0$, the following equation,  which
is an anisotropic generalization of (\ref{eq:general-eigen-prob}),
$$
-\nabla \cdot
\left(\begin{bmatrix}g(|\partial_{x}u_{0}|) & 0\\
0 & g(|\partial_{y}u_{0}|)
\end{bmatrix}\nabla u^{i}\right) = \lambda u^i,
$$
becomes exactly the Laplace equation in Grady's model which
is the Euler-Lagrange equation of the energy
\begin{equation}
\int_{\Omega}(\nabla u)^{T}\begin{bmatrix}g(|\partial_{x}u_{0}|) & 0\\
0 & g(|\partial_{y}u_{0}|)
\end{bmatrix}\nabla u .
\label{eq:bv-energy}
\end{equation}
While the proper definition of functional (\ref{eq:bv-energy})
for general $u$, $u_{0}$ possibly in $BV$ is missing, we can still
define it for a simpler case as in Proposition 5.2. Then, there is
no unique minimizer of this energy, and, as stated in Proposition
5.2, any minimizer of the above energy in $SBV$ yields the minimum
value 0 in the ideal case (with proper $g$ and $u_{0}$ as in Proposition
5.2). While the eigenvalue method considers the minimizer of the above
energy on the admissible set $\left\{ u:u|_{\partial\Omega}=0,\int u^{2}\ dx=1\right\} $,
the Random Walk method considers the minimizer satisfying boundary
conditions in (\ref{eq:Grady-model}). Both gives piecewise-constant
minimizers that can be used for image segmentation.

%
%
%
%
%

\section{Concluding remarks}
\label{SEC:conclusion}

We have introduced an eigenvalue problem of an anisotropic differential operator
as a tool for image segmentation. It is a continuous and anisotropic generalization 
of some commonly used, discrete spectral clustering models for image
segmentation. The continuous formulation of the eigenvalue problem allows
for accurate and efficient numerical implementation, which is crucial
in locating the boundaries in an image. An important observation
from numerical experiment is that non-trivial, almost piecewise constant
eigenfunctions associated with very small eigenvalues exist, and
this phenomenon seems to be an inherent property of the model.
These eigenfunctions can be used as the basis for image segmentation
and edge detection. The mathematical theory behind this is still unknown
and will be an interesting topic for future research.

We have implemented our model with a finite element method and shown
that anisotropic mesh adaptation is essential to the accuracy and efficiency
for the numerical solution of the model. Numerical tests on segmentation
of synthetic, natural or texture images based on computed eigenfunctions
have been conducted. It has been shown that the adaptive mesh implementation
of the model can lead to a significant gain in efficiency. 
Moreover, numerical results also show that the anisotropic nature of the
model can enhance some nontrivial regions of eigenfunctions which may
not be captured by a less anisotropic or an isotropic model.



\def\cprime{$'$}

\end{document}